\documentclass[12pt]{article}
\usepackage{amsthm}
\usepackage{amssymb}
\usepackage{mathrsfs}
\begin{document}
\renewcommand{\thefootnote}{\fnsymbol{footnote}}
\newfont{\bg}{cmr10 scaled\magstep4}
\newcommand{\bigzerou}{\smash{\lower1.7ex\hbox{\bg 0}}}
\newcommand{\bigzerol}{\smash{\hbox{\bg 0}}}
\newcommand{\bsquare}{\hbox{\rule{6pt}{6pt}}}
\newcommand{\isom}{\stackrel{\sim}{\to}}
\newcommand{\End}{\mbox{End}}
\newcommand{\Pic}{\mbox{Pic}}
\begin{center}
{\large\bf Standard models of abstract intersection theory}\\
{\large\bf for operators in Hilbert space}
\end{center}
\begin{center}
Grzegorz Banaszak$^{\ast}$\footnote[0]
{\noindent
\hspace{-0.7cm}
\noindent
\hspace{-0.7cm}
\noindent
\hspace{-0.7cm}
$^{\ast}$Supported by a research grant of the NCN (National Center for Science of Poland)}
and Yoichi Uetake
\end{center}
\begin{quotation}
\noindent
{\small
{\bf Abstract.} For an operator in a possibly infinite-dimensional Hilbert space of a certain class, we set down axioms of an abstract intersection theory, from which the Riemann hypothesis regarding the spectrum of that operator follows. 
In our previous paper [BU] we constructed a GNS (Gelfand-Naimark-Segal) model of abstract intersection theory.
In this paper we propose another model, which we call a standard model of abstract intersection theory. We show that there is a standard 
model of abstract intersection theory for a given operator if and only if the Riemann hypothesis and semi-simplicity hold for that operator. (For the definition
of semi-simplicity of an operator in Hilbert space, see the definition in Introduction.) 
We show this result under a condition for a given operator which is much weaker than the condition in the previous paper. The operator satisfying this condition
can be constructed by the method of automorphic scattering in [U].

Combining this with a result from [U],
we can show that an Dirichlet $L$-function, including the Riemann zeta-function, satisfies the Riemann hypothesis and its all nontrivial zeros are simple if 
and only if there is a corresponding standard model of abstract intersection theory. 
Similar results can be proven for GNS models since the same technique of
proof for standard models can be applied.\\
}
\end{quotation}
{\bf 1. Introduction}\\
\\
\indent
In the 1940s Weil [W1] developed an intersection theory on surfaces over finite fields to apply it to the proof of the Riemann hypothesis for curves over finite fields (one-variable function fields over finite fields).

In this paper we introduce axioms ((AIT1)--(AIT3) in $\S$3) of abstract intersection theory for an operator 
in a possibly infinite-dimensional Hilbert space, 
which are analogous to Weil's theory. We consider a collection ${\Bbb A}{\Bbb I}{\Bbb T}$ that consists of a vector space, its specific vectors and some maps, satisfying these axioms. 
From this collection one can derive the Riemann hypothesis regarding the spectrum of that operator. 
Therefore we call ${\Bbb A}{\Bbb I}{\Bbb T}$ an abstract intersection theory.  

Let $H$ be a possibly infinite-dimensional ${\Bbb C}$-Hilbert space.
Let $A\colon H\supset \mbox{dom}(A) \to H$ be a ${\Bbb C}$-linear operator acting on $H$. Here $\mbox{dom}(A)$ denotes
the domain of the operator $A$. We assume that its spectrum $\sigma(A)$ consists only of the point spectrum $\sigma_p(A)$. That is, $\sigma(A)$ consists only of eigenvalues of $A$. 

We say that the operator $A$ satisfies the Riemann hypothesis (RH, shortly) if
$$ \mbox{Re}(s_i)=\frac12~\mbox{for~all}~s_i\in \sigma(A)=\sigma_p(A). $$

We say that the operator $A$ is semi-simple if
$$ \nu(s_i)=1~\mbox{for~all}~s_i\in \sigma(A)=\sigma_p(A). $$
Here $\nu(s_i)$ is the Riesz index of $s_i$. For its definition see the paragraph preceding the conditions (OP1)--(OP5) in $\S$2, which $A$ is assumed to
satisfy. All these conditions are satisfied by an operator $A$ obtained from automorphic scattering theory [U], which gives a spectral interpretation
of a certain Dirichlet $L$-function, including the Riemann zeta-function. See Remark 2.1 (4) in $\S$2. 

In our previous work [BU], we showed ${\Bbb A}{\Bbb I}{\Bbb T} \Rightarrow$ RH. We also constructed a model
${\Bbb A}{\Bbb I}{\Bbb T}_{\mbox{\tiny GNS}}$ of abstract intersection theory based on an analogue of the GNS (Gelfand-Naimark-Segal) representation.
We call ${\Bbb A}{\Bbb I}{\Bbb T}_{\mbox{\tiny GNS}}$ a GNS model of abstract intersection theory.
We showed ${\Bbb A}{\Bbb I}{\Bbb T}_{\mbox{\tiny GNS}} \Leftrightarrow$ RH, assuming the semi-simplicity of $A$ ([BU, Theorem 3.1]).

We observe that there is some freedom in constructing models of abstract intersection theory to investigate the spectrum of operators in Hilbert space and
nontrivial zeros of corresponding Dirichlet $L$-functions.
In this paper we propose another new model ${\Bbb A}{\Bbb I}{\Bbb T}_m$, which we call a standard model of abstract intersection theory. 
This model is hinted by the K\"unneth formula for $\ell$-adic cohomology in the setting of the classical intersection theory.
For this model we show ${\Bbb A}{\Bbb I}{\Bbb T}_m \Leftrightarrow$ RH \& semi-simplicity (Theorem 5.2 (2)). 
The technique for proving this statement can also be applied
to ${\Bbb A}{\Bbb I}{\Bbb T}_{\mbox{\tiny GNS}}$ in the previous paper and one can show ${\Bbb A}{\Bbb I}{\Bbb T}_{\mbox{\tiny GNS}}
\Leftrightarrow$ RH \& semi-simplicity (Theorem 5.3). Therefore we significantly strengthen our previous results in [BU] for both GNS and standard models,
dropping the semi-simplicity assumption (the condition (OP3-b) in [BU]). The condition (OP3-b) in this paper is much weaker and is satisfied by  
operators coming from scattering theory for Dirichlet $L$-functions [U].

As a consequence of Theorems 5.2 and 5.3 combined with the results in [U] from automorphic scattering theory, we can show that an Dirichlet $L$-function, 
including the Riemann zeta-function, satisfies the RH and its all nontrivial zeros are simple if and only if there is a corresponding standard model 
${\Bbb A}{\Bbb I}{\Bbb T}_m$ (or GNS model ${\Bbb A}{\Bbb I}{\Bbb T}_{\mbox{\tiny GNS}}$) of abstract intersection theory (Theorem 5.4).

The plan of this paper is as follows.

In $\S$2 we define an analogue of the classical Frobenius morphism for the operator $A$. The spectrum of 
this analogue is similar to that of the classical Frobenius morphism if the operator $A$ satisfies the Riemann hypothesis. The introduction of this analogue
is also hinted by Weil's explicit formulas [W2].

In $\S$3 we introduce a general notion of abstract intersection theory ${\Bbb A}{\Bbb I}{\Bbb T}$ and set down its axioms ((AIT1), (AIT2) and (AIT3)).

In $\S$4 we construct a specific example of abstract intersection theory, which we call a standard model ${\Bbb A}{\Bbb I}{\Bbb T}_m$, using analogy with the classical K\"unneth formula for $\ell$-adic cohomology. 

In $\S$5 we state our main theorems (Theorems 5.2, 5.3 and 5.4).

In $\S$6 we show that there is a strong analogy between Weil's approach to zeta-functions for curves over finite fields and our approach to Dirichlet $L$-functions.
For Weil's intersection theory, see also Grothendieck [Gro], Monsky [Mon] and Serre [S].

We should note that there is a program by Connes and Marcolli (and Consani) [CM] to adapt Weil's proof of RH for function fields to the case of number fields.
See also Connes [C].
There is also a conjectural cohomological approach by Deninger [D] toward the interpretation of $L$-functions analogous 
to the etale cohomology theory of varieties over finite fields.\\ 
\\ 
\noindent
{\bf 2. An analogue of the Frobenius morphism for the operator $A$}\\
\\
\indent
Let $H$ be a possibly infinite-dimensional ${\Bbb C}$-Hilbert space.
If $H$ is infinite-dimensional we assume that $H$ is separable. 
Let $A\colon H\supset \mbox{dom}(A) \to H$ be a possibly unbounded operator on $H$. 

If $s_i \in \sigma(A)$ is an isolated spectrum point, one can take a small enough bounded domain $\Delta$ 
of ${\Bbb C}$ such that $\{ s_i \} \Subset \Delta$ (i.e.,\,$\{ s_i \} \subset \Delta^{\circ}$) and $\overline{\Delta}\cap (\sigma(A)-\{ s_i \})=\emptyset$. 
Then one can define the Riesz projection $P_{\{ s_i \}}\colon H\to H$ by
$$ P_{\{ s_i \}}:=\frac{1}{2\pi i}\oint_{\partial \Delta} (sI-A)^{-1}ds. $$
Here $I\colon H \to H$ is the identity operator on $H$. $P_{\{s_i\}}$ is a bounded operator on $H$.

For $s_i \in \sigma_p(A)$, the Riesz index $\nu(s_i)$ of $s_i$ is defined as the smallest positive number $\leq \infty$ such that 
$$ \mbox{Ker}((s_i I-A)^{\nu(s_i)}) = \mbox{Image}(P_{\{ s_i \}}). $$

Let $\mbox{mult}(s_i)=\dim_{\Bbb C} \mbox{Image}(P_{\{ s_i \}})$, which we call the (algebraic) multiplicity of $s_i \in \sigma_p(A)$.\\

We assume the following properties of $A$.\\
\\
{\bf (OP1)} $A$ is closed.\\
\\
{\bf (OP2)} The spectrum $\sigma(A)$ consists only of the point spectrum (i.e.,\,eigenvalues) $\sigma_p(A)$, i.e.,\,$\sigma(A)=\sigma_p(A)$, 
which accumulates at most at infinity.\\
\\
{\bf (OP3)} (a) $\mbox{Image}(P_{\{s_i\}})$ is finite-dimensional for any $s_i \in \sigma_p(A)$.\\
\indent \hspace{0.6cm} 
(b) $\nu(s_i)=\mbox{mult}(s_i)$ for any $s_i \in \sigma_p(A)$.\\ 
\\
{\bf (OP4)} $\sigma(A)\subset \Omega_{\infty}$, where 
$\Omega_{\infty}:=\{s\in {\Bbb C}\vert 0<\mbox{Re}(s)<1 \}$.\\
\\
{\bf (OP5)} 
$\mbox{Re}(s_i)<\frac12$ for some $s_i \in \sigma(A)$ if and only if there is $s_j \in \sigma(A)$ such that $\mbox{Re}(s_j)>\frac12$.\\
\\
{\it Remark 2.1.}\\
(1) (OP1) is needed when one applies Lemma 2.1 of [BU] on spectral decomposition.
Lemma 2.1 of [BU] is taken from Gohberg, Goldberg and Kaashoek [GoGK, XV.2, Theorem 2.1, p.\,326].\\
(2) In the previous paper [BU], the condition (OP3-b) was the simi-simplicity $\nu(s_i)=1$ ($s_i \in \sigma(A)$). The above stated (OP3-b) is a much weaker condition. This condition says that each spectrum (eigenvalue) of $A$ has just one corresponding Jordan block.
Actually this is satisfied in the construction using automorphic scattering theory [U]. See Remark 2.1\,(4) below.\\
(3) The above (OP5) is (OP5-a) in [BU]. (OP5-b) in [BU], which is necessary for the construction of GNS models of abstract intersection theory is not
necessary for the construction of standard models in this paper. (OP5-b) in [BU] is used to keep the space $V$ an ${\Bbb R}$-linear space in the GNS model.
In the standard model we apply the complexification $V_{\Bbb C}$ of $V$ instead (see $\S$3).
(OP5-b) in [BU] is satisfied by an operator $A$ constructed in [U] (see Remark 2.1\,(4) below). \\
(4) Let $\Gamma$ be a congruence subgroup of $SL_2({\Bbb Z})$ such that $\Gamma\setminus {\Bbb H}~(\simeq \Gamma\setminus SL_2({\Bbb R})/SO(2))$ is noncompact and has one cusp at $i\infty$. Here ${\Bbb H}$ denotes the upper half-plane. 
In [U] the second author constructed a scattering theory for automorphic forms on $\Gamma \setminus {\Bbb H}$. 
Furthermore he constructed an operator $A$ satisfying (OP1)--(OP5) whose (point) spectrum coincides with the nontrovial zeros
of the Dirichlet $L$-function $L(s,\chi)$ associated to $\Gamma$, counted with multiplicity:\,$s_i \in \sigma(A)(=\sigma_p(A))$ and $\nu(s_i)(=\mbox{mult}(s_i))=m_i$ 
if and only if $s_i$ is a nontrivial zero of $L(s,\chi)$ of order $m_i$. We call $s_i \in {\Bbb C}$ a nontrivial zero of the Dirichlet $L$-function $L(s,\chi)$ if $L(s_i,\chi)=0$ and $0 < \mbox{Re}(s_i) < 1$. See Theorem 4.1 (i) ($\Rightarrow$ (OP1)), (iii-a), (iii-b), (iii-c) ($\Rightarrow$ (OP2) and (OP3)) 
and (iv) ($\Rightarrow$ (OP4) and (OP5), and (OP5-b) in [BU]) of [U, p.\,455]. The theory of automorphic scattering was initiated by Pavlov-Faddeev [PavF] and 
then Lax-Phillips [LP] hinted by Gelfand [Ge]. \hfill $\Box$\\
\\
\indent
Now for $Y>0$ let 
$$ \sigma_Y(A) := \{s\in \sigma(A)~\vert~\vert\mbox{Im}(s)\vert < Y\}. $$
Note that $\sigma_Y(A)$ is a finite set by (OP2).
Let the parameter space ${\cal Y}$ be defined by
$$ {\cal Y} := \{ Y > 0~\vert~\sigma_Y(A) \neq \emptyset \} - \{ \vert\mbox{Im}(s)\vert~\vert~s \in \sigma(A) \}. $$
Fix a function
$$ q\colon {\cal Y} \to (0,1) \cup (1,\infty). $$

Let $B(X)$ denote the set of bounded operators on a ${\Bbb C}$-Hilbert space $X$. 
By definition $T\colon X \supset \mbox{dom}(T) \to X$ is a bounded operator if $\mbox{dom}(T)=X$ and the operator norm $\Vert T \Vert < \infty$.    

Let $\Sigma_H$ be the set of closed subspaces of $H$. We will construct maps
$$ F_A\colon {\cal Y} \to B(H) $$
and 
$$ {\cal H}\colon {\cal Y} \to \Sigma_H $$
such that $F_A(Y)\colon H \to H$ satisfies the following conditions for each $Y \in {\cal Y}$.
$$ F_A(Y){\cal H}(Y)\subset {\cal H}(Y) \leqno{\mbox{{\bf (FROB-a)}}} $$
(i.e.,\,the subspace ${\cal H}(Y)$ is invariant for $F_A(Y)$). 

$$ \sigma(F_A(Y)\vert_{{\cal H}(Y)}) = \sigma_p(F_A(Y)\vert_{{\cal H}(Y)}) = \{ q(Y)^s \vert s \in \sigma_Y(A) \} \leqno{\mbox{{\bf (FROB-b)}}} $$
and
$$ \sigma(F_A(Y)) = \sigma_p(F_A(Y)) = \sigma(F_A(Y)\vert_{{\cal H}(Y)}) \cup \{0\}. $$
Note that $\sigma(F_A(Y)\vert_{{\cal H}(Y)})$ is a finite set since $\sigma_Y(A)$ is. 

The operator $F_A(Y)$ $(Y \in {\cal Y})$ is considered to be an analogue of the classical Frobenius morphism, since the spectrum of 
this analogue is similar to that of the classical Frobenius morphism if the operator $A$ satisfies the Riemann hypothesis.
It is also hinted by the spectral side of Weil's explicit formulas [W2] (see $\S$6).\\
\\
{\it Models} $F_{A,m}$ {\it and} ${\cal H}_m$ {\it of} $F_A$ {\it and} ${\cal H}$:\\
Now we construct the models $F_{A,m}\colon {\cal Y} \to B(H)$ and ${\cal H}_m\colon {\cal Y} \to \Sigma_H$ which satisfy (Frob-a) and (Frob-b). 
These models will constitute parts of a standard model ${\Bbb A}{\Bbb I}{\Bbb T}_m$ constructed in $\S$4. Let
$$ \Omega_Y := \{s\in {\Bbb C}\vert 0< \mbox{Re}(s)< 1,\,\vert \mbox{Im}(s) \vert < Y \} $$
for $Y \in {\cal Y}$. 
Note that $\Omega_Y \cap \sigma(A)=\sigma_Y(A)$ for $Y \in {\cal Y}$ by (OP4).
Note that for $Y \in {\cal Y}$,
$$ \sigma_Y(A) \Subset \Omega_Y~~(\mbox{i.e.,}\,\overline{\sigma_Y(A)}=\sigma_Y(A) \subset \Omega_Y^{\circ}=\Omega_Y)~~ \mbox{and}~~ \overline{\Omega_Y}\cap (\sigma(A)-\sigma_Y(A))=\emptyset. $$
Therefore, for $Y \in {\cal Y}$, the Riesz projection 
$$ P_{\sigma_Y(A)}\colon H \to H $$
can be well-defined by
$$ P_{\sigma_Y(A)}:=\frac{1}{2\pi i}\oint_{\partial \Omega_Y} (sI-A)^{-1}ds. $$
$P_{\sigma_Y(A)}$ is a bounded operator on $H$. Let ${\cal H}_m\colon {\cal Y} \to \Sigma_H$ be defined by 
$$ {\cal H}_m(Y):=\mbox{Image}(P_{\sigma_Y(A)}). $$
By (OP2) and (OP3-a), ${\cal H}_m(Y)$ is finite-dimensional for each $Y \in {\cal Y}$. 

Given $Y \in {\cal Y}$, let 
$$ F_{A,m}(Y)\colon H\supset \mbox{dom}(F_{A,m}(Y)) \to H $$ 
be defined by
$$ F_{A,m}(Y) x := \frac{1}{2\pi i}\Bigl(\oint_{\partial\Omega_Y} q(Y)^s (sI-A)^{-1} ds\Bigr)x $$
for
$$ x \in \mbox{dom}(F_{A,m}(Y)) := \{ x \in H \vert F_{A,m}(Y)x~\mbox{exists~in}~H \}. $$

Note that $\sigma_Y(A)$ is a bounded set. Thus, by Lemma 2.1 of [BU], ${\cal H}_m(Y) \subset \mbox{dom}(A)$ and $A{\cal H}_m(Y) \subset {\cal H}_m(Y)$.
Let $A\vert_{{\cal H}_m(Y)}\colon {\cal H}_m(Y) \to {\cal H}_m(Y) \subset H$ be the restriction of $A$ to ${\cal H}_m(Y)$. 
Since ${\cal H}_m(Y)$ is finite-dimensional, $A\vert_{{\cal H}_m(Y)}$ is a bounded operator, i.e.,\,$A\vert_{{\cal H}_m(Y)} \in B({\cal H}_m(Y))$.

Similarly, let $H(s_i):=\mbox{Image}(P_{\{ s_i \}})$. By (OP2) and (OP3-a), $H(s_i)$ is finite-dimensional. Again by Lemma 2.1 of [BU], $H(s_i) \subset 
\mbox{dom}(A)$ and $AH(s_i) \subset H(s_i)$. Let $A\vert_{H(s_i)}\colon H(s_i) \to H(s_i) \subset H$ be the restriction of $A$ to $H(s_i)$. 
By the same argument for $A\vert_{{\cal H}_m(Y)}$, we have $A\vert_{H(s_i)} \in B(H(s_i))$.\\
\\
{\bf Lemma 2.1.} {\it Suppose that $A$ satisfies} (OP1), (OP2) {\it and} (OP3-a). {\it Then}\\
(i) {\it For each $Y \in {\cal Y}$,} 
$$ \mbox{dom}(F_{A,m}(Y))=H. $$
{\it Furthermore, $F_{A,m}(Y)$ is a bounded operator on $H$, i.e.,\,$F_{A,m}(Y) \in B(H)$.}\\
(ii) {\it The subspace ${\cal H}_m(Y)$ is $F_{A,m}(Y)$-invariant}: 
$$ F_{A,m}(Y){\cal H}_m(Y)\subset {\cal H}_m(Y). $$
{\it That is, $F_{A,m}$ satisfies} (Frob-a).\\
(iii) {\it For each $Y \in {\cal Y}$, we have} 
$$ \sigma(F_{A,m}(Y)\vert_{{\cal H}_m(Y)}) = \sigma_p(F_{A,m}(Y)\vert_{{\cal H}_m(Y)}) = \{ q(Y)^s \vert s \in \sigma_Y(A) \} $$
{\it and}
$$ \sigma(F_{A,m}(Y)) = \sigma_p(F_{A,m}(Y)) = \sigma(F_{A,m}(Y)\vert_{{\cal H}_m(Y)}) \cup \{0\}. $$
{\it That is, $F_{A,m}$ satisfies} (Frob-b).\\
(iv) {\it Let $t\colon {\cal Y} \to {\Bbb R}-\{0\}$ be defined by $t(Y):=\log q(Y)$ $($i.e.,\,$e^{t(Y)}=q(Y)$$)$ for 
$Y \in {\cal Y}$. 
For each $Y \in {\cal Y}$, we have}
$$ F_{A,m}(Y) = e^{t(Y)A\vert_{{\cal H}_m(Y)}}P_{\sigma_Y(A)}=\sum_{n=0}^{\infty} \frac{t(Y)^n}{n!}A\vert_{{\cal H}_m(Y)}^n P_{\sigma_Y(A)} 
=\sum_{s_i \in \sigma_Y(A)} e^{t(Y)A\vert_{H(s_i)}} P_{\{s_i \}}. $$
(v) {\it Suppose further that $A$ satisfies} (OP3-b).
{\it Then, with respect to an appropriate basis of $H(s_i)$, $e^{t(Y)A\vert_{H(s_i)}}$ is written as}
$$ e^{t(Y)A\vert_{H(s_i)}} = N(s_i) $$
{\it with}      
$$ N(s_i) = 
\left(\begin{array}{ccccc}
\frac{t(Y)^0e^{t(Y)s_i}}{0!} & \frac{t(Y)^1e^{t(Y)s_i}}{1!} & \cdots                       & \cdots                       & \frac{t(Y)^{m_{i}-1}e^{t(Y)s_i}}{(m_{i}-1)!} \\ 
                             & \frac{t(Y)^0e^{t(Y)s_i}}{0!} & \frac{t(Y)^1e^{t(Y)s_i}}{1!} & \cdots                       & 
\frac{t(Y)^{m_{i}-2}e^{t(Y)s_i}}{(m_{i}-2)!} \\
                             &                              & \ddots                       & \ddots                       & \vdots                       \\
                             &                              &                              & \frac{t(Y)^0e^{t(Y)s_i}}{0!} & \frac{t(Y)^1e^{t(Y)s_i}}{1!} \\ 
\bigzerol                    &                              &                              &                              & \frac{t(Y)^0e^{t(Y)s_i}}{0!}
\end{array}\right) \in M_{m_i}({\Bbb C}). $$
{\it Here} $m_{i}=\nu(s_i)(=\mbox{mult}(s_i))$.\\
\\
\noindent
{\it Proof.}
Let $K(Y)=\mbox{Ker}(P_{\sigma_Y(A)})$. Then by Lemma 2.1 of [BU], $K(Y)$ is $A$-invariant in the sense that $A(K(Y)\cap\mbox{dom}(A))\subset K(Y)$.
Thus one can define $A\vert_{K(Y)}\colon K(Y) \supset \mbox{dom}(A\vert_{K(Y)}) \to K(Y)$, the restriction of $A$ to $K(Y)$. Then we have $\sigma(A\vert_{K(Y)})=\sigma(A)-\sigma_Y(A)$ (Lemma 2.1 [BU]).
We also have
$$ A=\left(\begin{array}{cc}
A\vert_{{\cal H}_m(Y)} & 0 \\ 
0 & A\vert_{K(Y)}  
\end{array}\right) $$
on $H={\cal H}_m(Y) \oplus K(Y)$. Note that the direct sum $\oplus$ does not necessarily mean the orthogonal sum.

By (OP2) and (OP3-a), $(sI-A)^{-1}$ is meromorphic in the whole ${\Bbb C}$-plane.
However, since $(sI-A\vert_{K(Y)})^{-1}$ is holomorphic in $\Omega_Y$, we have by the functional calculus for the bounded operator $A\vert_{{\cal H}_m(Y)}$
\begin{eqnarray*}
F_{A,m}(Y) &=& \frac{1}{2\pi i}\oint_{\partial\Omega_Y} q(Y)^s (sI-A)^{-1} ds \\ 
  &=& \frac{1}{2\pi i}\oint_{\partial\Omega_Y} e^{t(Y)s} (sI-A)^{-1} ds \\
  &=& \frac{1}{2\pi i}\oint_{\partial\Omega_Y} e^{t(Y)s} \Biggl(sI-
       \left(\begin{array}{cc}
       A\vert_{{\cal H}_m(Y)} & 0 \\ 
       0 & A\vert_{K(Y)}  
       \end{array}\right)\Biggr)^{-1} ds \\
  &=& \frac{1}{2\pi i}\oint_{\partial\Omega_Y} e^{t(Y)s}
       \left(\begin{array}{cc}
       (sI-A\vert_{{\cal H}_m(Y)})^{-1} & 0 \\ 
       0 & (sI-A\vert_{K(Y)})^{-1}  
       \end{array}\right) ds \\
  &=& \frac{1}{2\pi i}\oint_{\partial\Omega_Y} e^{t(Y)s} (sI-A\vert_{{\cal H}_m(Y)})^{-1} P_{\sigma_Y(A)} ds\\
  &=& e^{t(Y)A\vert_{{\cal H}_m(Y)}} P_{\sigma_Y(A)}\\
  &=& \left(\begin{array}{cc}
       e^{t(Y)A\vert_{{\cal H}_m(Y)}} & 0 \\ 
       0 & 0  
      \end{array}\right),
\end{eqnarray*}
which shows (i) and (ii). By Lemma 2.1 of [BU], we have $\sigma(A\vert_{{\cal H}_m(Y)})=\sigma_Y(A)$.
Applying the spectral mapping theorem to the bounded operator $A\vert_{{\cal H}_m(Y)}$ (recall that $\dim_{\Bbb C} {\cal H}_m(Y)<\infty$), this also shows (iii).

Note that 
$$ P_{\sigma_Y(A)}=\bigoplus_{s_i \in \sigma_Y(A)} P_{\{s_i\}} $$
and
$$ {\cal H}_m(Y) = \bigoplus_{s_i \in \sigma_Y(A)} H(s_i). $$
Here $\bigoplus$ denotes the (not necessarily orthogonal) direct sum.
Therefore we have
$$ F_{A,m}(Y) = \frac{1}{2\pi i}\sum_{s_i \in \sigma_Y(A)} \oint_{\partial\Omega_Y} q(Y)^s (sI-A\vert_{H(s_i)})^{-1} P_{\{s_i\}} ds
= \sum_{s_i \in \sigma_Y(A)} e^{t(Y)A\vert_{H(s_i)}} P_{\{s_i \}}. $$
From this (iv) follows. 

Note that by Lemma 2.1 of [BU] we have $\sigma(A\vert_{H(s_i)})=\{ s_i \}$. Thus, by (OP3-b), $A\vert_{H(s_i)}$ is written 
with respect to an appropriate basis of $H(s_i)$ as
$$ A\vert_{H(s_i)} = M(s_i) $$
with      
$$ M(s_i) =
\left(\begin{array}{ccccc}
s_i       & 1   &        &        & \bigzerou \\ 
          & s_i & 1      &        &           \\
          &     & \ddots & \ddots &           \\
          &     &        & s_i    & 1         \\ 
\bigzerol &     &        &        & s_i 
      \end{array}\right) \in M_{m_i}({\Bbb C}). $$
Here $m_{i}=\nu(s_i)$.

Note that
$$ (sI-M(s_i))^{-1}= 
\left(\begin{array}{ccccc}
\frac{1}{s-s_i} & \frac{1}{(s-s_i)^2} & \cdots              & \cdots          & \frac{1}{(s-s_i)^{m_{i}}}   \\ 
                & \frac{1}{s-s_i}     & \frac{1}{(s-s_i)^2} & \cdots          & \frac{1}{(s-s_i)^{m_{i}-1}} \\
                &                     & \ddots              & \ddots          & \vdots                    \\
                &                     &                     & \frac{1}{s-s_i} & \frac{1}{(s-s_i)^2}       \\ 
\bigzerol       &                     &                     &                 & \frac{1}{s-s_i} 
\end{array}\right). $$
Note that $q(Y)^s=e^{t(Y)s}=\sum_{n=0}^{\infty}\frac{t(Y)^ne^{t(Y)s_i}}{n!}(s-s_i)^n$.
From this (v) follows by using the residue theorem. \hfill $\Box$\\
\\
\noindent
{\bf 3. Abstract intersection theory and its axioms}\\
\\
\indent
Let $V$ be an ${\Bbb R}$-linear space endowed with a symmetric ${\Bbb R}$-bilinear form $\beta\colon V\times V \to {\Bbb R}$.
Denote by $V_{\Bbb C}$ the complexification of $V$ given by $V_{\Bbb C}=V \otimes_{\Bbb R} {\Bbb C}$.
To simplify the notation, we identify $v \otimes \alpha$ with $\alpha v$ for $v \in V$ and $\alpha \in {\Bbb C}$. Therefore we have $V \subset V_{\Bbb C}$.  
Then one can define the complexification $\beta_{\Bbb C}\colon V_{\Bbb C} \times V_{\Bbb C} \to {\Bbb C}$ of $\beta$ by
$$ \beta_{\Bbb C}(\alpha_1 v_1, \alpha_2 v_2):=\alpha_1 \overline{\alpha_2} \beta(v_1, v_2) \quad (v_1, v_2 \in V, \alpha_1, \alpha_2 \in {\Bbb C}). $$
It is easy to check that $\beta_{\Bbb C}(\alpha w_1, w_2)=\alpha \cdot \beta_{\Bbb C}(w_1, w_2)$ and 
$\beta_{\Bbb C}(w_2, w_1)=\overline{\beta_{\Bbb C}(w_1, w_2)}$ 
for $w_1, w_2 \in V_{\Bbb C}$ and $\alpha \in {\Bbb C}$.

Let $\mbox{End}_{\Bbb C}(V_{\Bbb C})$ denote the set of ${\Bbb C}$-linear operators $T\colon V_{\Bbb C} \supset \mbox{dom}(T) \to V_{\Bbb C}$
such that $\mbox{dom}(T)=V_{\Bbb C}$.
Suppose that there are nonzero vectors $v_{01}$, $v_{10}$ and $h_a$ in $V$, 
maps $v_{\delta}\colon {\cal Y} \to V_{\Bbb C}$ and $\Phi_A \colon {\cal Y} \to\mbox{End}_{\Bbb C}(V_{\Bbb C})$
which satisfy the conditions listed below ((AIT1)--(AIT3)). We call a collection 
$$ {\Bbb A}{\Bbb I}{\Bbb T}=(V,v_{01},v_{10},v_{\delta},h_a,\beta,\Phi_A,F_A, {\cal H}) $$
which satisfies these conditions an {\it abstract intersection theory}. 
The map $\Phi_A$ is associated with the operator $A$ in $\S$2. 
$F_A\colon {\cal Y} \to B(H)$ along with ${\cal H}\colon {\cal Y} \to \Sigma_H$ is an analogue of the Frobenius morphism defined in $\S$2, which satisfies 
(Frob-a) and (Frob-b). $F_A$ is related with $\Phi_A$ by the axiom (AIT3).\\
\\ 
\noindent
{\bf (AIT1)} (a) $\beta(y,x)=\beta(x,y) \in {\Bbb R}$~ for~ $x, y \in V$. 
$\beta_{\Bbb C}(y,x)=\overline{\beta_{\Bbb C}(x,y)} \in {\Bbb C}$~ for~ $x, y \in V_{\Bbb C}$.\\ 
\indent \hspace{0.8cm}
(b) $\beta(v_{01},v_{01})=0$.\quad 
(c) $\beta(v_{10},v_{10})=0$.\quad 
(d) $\beta(v_{01},v_{10})=1$.\\
\indent \hspace{1.5cm}
For each $Y \in {\cal Y}$ and all $n\geq 0$:\\
\indent \hspace{0.8cm}
(e) $\beta_{\Bbb C}(\Phi_A(Y)^n v_{\delta}(Y), v_{01})=1$.\quad  
(f) $\beta_{\Bbb C}(\Phi_A(Y)^n v_{\delta}(Y), v_{10})=O(q(Y)^n)$.\\
\indent \hspace{0.8cm}
(g) $\beta_{\Bbb C}(\Phi_A(Y)^n v_{\delta}(Y), \Phi_A(Y)^n v_{\delta}(Y))=O(q(Y)^n)$.\\
\\
{\bf (AIT2)} For $x \in V$, if $\beta(x,h_a)=0$ then $\beta(x,x)\leq 0$.\\
\\
\noindent
Note that (AIT1-e)--(AIT1-g) are assumed to hold for each $Y \in {\cal Y}$.
The Bachmann-Landau notation $O(q(Y)^n)$ in (AIT1) is with respect to $n\gg 0$ for $q(Y)$ with $Y\in {\cal Y}$ fixed.
We call (AIT2) the {\it Hodge property}, and $h_a$ a {\it Hodge vector}.\\    
\\
\noindent
{\bf Lemma 3.1.} 
{\it Under the assumptions} (AIT1-a)--(AIT1-d) {\it and} (AIT2), 
{\it we have}
$$ \beta(x,x) \leq 2 \beta(x,v_{01})\beta(x,v_{10})\quad
(x \in V). $$
{\it Proof.} See the proof of [BU, Lemma 3.1]. \hfill $\Box$\\
\\
\indent
Let the ${\Bbb R}$-bilinear form $\langle \cdot, \cdot \rangle_V\colon V \times V \to {\Bbb R}$ be defined by 
$$ \langle x, y \rangle_{V}:=\beta(x,v_{01})\beta(v_{10},y)+\beta(x,v_{10})\beta(v_{01},y)-\beta(x,y) \leqno{{\bf (\ast)}} $$
for $x, y \in V$.
By Lemma 3.1, $\langle \cdot, \cdot \rangle_{V}$ is positive semidefinite, 
i.e.,\,$\langle x,x \rangle_{V} \geq 0$ for $x\in V$. 
Indeed, as we will see soon below ((IP-b), (IP-c)), this bilinear form must be positive
{\it semi}definite, not positive definite.

We obtain the complexification $\langle \cdot, \cdot \rangle_{V_{\Bbb C}}\colon V_{\Bbb C} \times V_{\Bbb C} \to {\Bbb C}$ of 
$\langle \cdot, \cdot \rangle_V\colon V \times V \to {\Bbb R}$ by 
$$ \langle \alpha_1 v_1, \alpha_2 v_2 \rangle_{V_{\Bbb C}}:=\alpha_1 \overline{\alpha_2} \langle v_1, v_2 \rangle_V $$ 
for $v_1, v_2 \in V$ and $\alpha_1, \alpha_2 \in {\Bbb C}$.\\
\\
{\bf Lemma 3.2.} {\it $\langle \cdot, \cdot \rangle_{V_{\Bbb C}}$ is positive semidefinite, i.e.,\,$\langle x, x \rangle_{V_{\Bbb C}} \geq 0$ 
for all $x \in V_{\Bbb C}$.}\\
\\
{\it Proof.} Since for $x, y \in V$ and $t \in {\Bbb R}$,
$$ \langle tx+y, tx+y \rangle_V = \langle x, x \rangle_V t^2 + 2 \langle x, y \rangle_V t + \langle y, y \rangle_V \geq 0, $$
we have the Cauchy-Schwarz inequality for $\langle \cdot, \cdot \rangle_V$ 
$$ \vert \langle x,y \rangle_V \vert \leq \sqrt{\langle x,x \rangle_V \langle y,y \rangle_V}\qquad (x,y\in V), $$
provided that $\langle x, x \rangle_V \neq 0$. If $\langle x, x \rangle_V=0$ then $\langle x, y \rangle_V$ also must be zero. Therefore we have the
Cauchy-Schwarz inequality for $\langle \cdot, \cdot \rangle_V$ for any $x, y \in V$. 

Let ${\cal V}$ be a basis of $V$. Split ${\cal V}$ into two disjoint sets 
${\cal V}=\{ u_i \}_{i \in I} \cup \{ v_j \}_{j \in J}$ so that $\langle u_i, u_i \rangle_V=0$ and 
$\langle v_j, v_j \rangle_V\neq 0$. Note that ${\cal V}$ is also a basis of $V_{\Bbb C}$ with the same properties that $\langle u_i, u_i \rangle_{V_{\Bbb C}}=0$ and 
$\langle v_j, v_j \rangle_{V_{\Bbb C}}\neq 0$. Therefore any $x \in V_{\Bbb C}$ can be written as
$$ x = \sum_{i \in I_x} \alpha_{x,i} u_i + \sum_{j \in J_x} \alpha_{x,j} v_j $$
for some {\it finite} subsets $I_x \subset I$ and $J_x \subset J$ with $\alpha_{x,i}, \alpha_{x,j} \in {\Bbb C}$. 

Apply the Gram-Schmidt process to $\{ v_j \}_{j \in J_x}$ in $V$ to obtain an orthonormal set
$\{ e_j \}_{j \in J_x}$ in $V$. Then we have
$$ x = \sum_{i \in I_x} \alpha_{x,i} u_i + \sum_{j \in J_x} \alpha_{x,j}^{\prime} e_j $$
for some $\alpha_{x,j}^{\prime} \in {\Bbb C}$.

From the Cauchy-Schwarz inequality for $\langle \cdot, \cdot \rangle_V$, 
we have $\langle u_{i_1}, u_{i_2} \rangle_V=\langle u_{i_1}, u_{i_2} \rangle_{V_{\Bbb C}}=0$ for $i_1, i_2 \in I_x$ and
$\langle u_i, e_j \rangle_V=\langle u_i, e_j \rangle_{V_{\Bbb C}}=0$
for $i \in I_x$ and $j \in J_x$. Thus it is easy to see that $\langle x, x \rangle_{V_{\Bbb C}} \geq 0$. \hfill $\Box$\\
\\
\indent
Note that $\langle \cdot, \cdot \rangle_{V_{\Bbb C}}$ is compatible with $\beta_{\Bbb C}$, i.e.,\,we have
$$ \langle x, y \rangle_{V_{\Bbb C}}=\beta_{\Bbb C}(x,v_{01})\beta_{\Bbb C}(v_{10},y)+\beta_{\Bbb C}(x,v_{10})\beta_{\Bbb C}(v_{01},y)
-\beta_{\Bbb C}(x,y). \leqno{{\bf (\ast\ast)}} $$

It is easy to see that from (AIT1), ($\ast$) and ($\ast\ast$) the following conditions follow for any $Y \in {\cal Y}$.\\
\\
{\bf (IP)} (a) $\langle y, x \rangle_V=\langle x, y \rangle_V \in {\Bbb R}$ for $x, y \in V$.
$\langle y, x \rangle_{V_{\Bbb C}}=\overline{\langle x, y \rangle_{V_{\Bbb C}}} \in {\Bbb C}$ for $x, y \in V_{\Bbb C}$.\\
\indent \hspace{0.2cm}
(b) $\langle v_{01},v_{01} \rangle_V=0$.\quad 
(c) $\langle v_{10},v_{10} \rangle_V=0$.\quad 
(d) $\langle v_{01},v_{10} \rangle_V=0$.\\
\indent \hspace{1.0cm}
For each $Y \in {\cal Y}$ and all $n \geq 0$:\\
\indent \hspace{0.2cm}  
(e) $\langle \Phi_A(Y)^n v_{\delta}(Y), v_{01} \rangle_{V_{\Bbb C}}=0$.\quad 
(f) $\langle \Phi_A(Y)^n v_{\delta}(Y), v_{10} \rangle_{V_{\Bbb C}}=0$.\\
\indent \hspace{0.2cm}
(g) $\langle \Phi_A(Y)^n v_{\delta}(Y), \Phi_A(Y)^n v_{\delta}(Y) \rangle_{V_{\Bbb C}}=O(q(Y)^n)$.\\ 
\\
\noindent
{\bf Lemma 3.3.} {\it For $\langle \cdot, \cdot \rangle_{V_{\Bbb C}}$, we have the Cauchy-Schwarz inequality}: 
$$ \vert \langle x,y \rangle_{V_{\Bbb C}} \vert \leq 
\sqrt{\langle x,x \rangle_{V_{\Bbb C}} \langle y,y \rangle_{V_{\Bbb C}}}\qquad (x,y\in V_{\Bbb C}). $$
\\
{\it Proof.} Let $\lambda=\langle x, x \rangle_{V_{\Bbb C}}$. By Lemma 3.2 we have $\lambda \geq 0$. Note that (e.g.,\,MacCluer [Mac, Exercise 1.7, p.\,24])
$$ 0 \leq \langle \lambda y - \langle y, x \rangle_{V_{\Bbb C}} x, \lambda y - \langle y, x \rangle_{V_{\Bbb C}} x \rangle_{V_{\Bbb C}}
= \lambda \{ \lambda \langle y, y \rangle_{V_{\Bbb C}} - \vert \langle x, y \rangle_{V_{\Bbb C}} \vert^2 \}. $$ 
Therefore if $\lambda > 0$ we have the inequality. Suppose $\lambda=0$. For the basis ${\cal V}$ in the proof of Lemma 3.2,
$$ x = \sum_{i \in I_x} \alpha_{x,i} u_i + \sum_{j \in J_x} \alpha_{x,j} v_j $$
for some {\it finite} subsets $I_x \subset I$ and $J_x \subset J$ with $\alpha_{x,i}, \alpha_{x,j} \in {\Bbb C}$. Applying the Gram-Schmidt process to 
$\{ v_j \}_{j \in J_x}$ in $V$, obtain an orthonormal set $\{ e_j \}_{j \in J_x}$ in $V$. Then as in the proof of Lemma 3.2 we have for some $\alpha^{\prime}_{x,j}$
$$ x = \sum_{i \in I_x} \alpha_{x,i} u_i + \sum_{j \in J_x} \alpha^{\prime}_{x,j} e_j. $$
Since $\lambda=0$ we have $\alpha^{\prime}_{x,j}=0$.  Therefore
$$ x = \sum_{i \in I_x} \alpha_{x,i} u_i. $$
Similarly, $y$ can be expressed as
$$ y = \sum_{i \in I_y} \alpha_{y,i} u_i + \sum_{j \in J_y} \alpha_{y,j} v_j $$
for some {\it finite} subsets $I_y \subset I$ and $J_y \subset J$ with $\alpha_{y,i}, \alpha_{y,j} \in {\Bbb C}$.
Since $\langle u_i, u_i \rangle_V=0$ for $i \in I_x$, we have, by the Cauchy-Schwarz inequality for $\langle \cdot, \cdot \rangle_V$, 
$\langle u_{i_1}, u_{i_2} \rangle_V=\langle u_{i_1}, u_{i_2} \rangle_{V_{\Bbb C}}=0$ for $i_1 \in I_x$ and $i_2 \in I_y$ 
and $\langle u_{i}, v_{j} \rangle_V=\langle u_{i}, v_{j} \rangle_{V_{\Bbb C}}=0$ for $i \in I_x$ and $j \in J_y$. 
Thus we have $\langle x, y \rangle_{V_{\Bbb C}}=0$. \hfill $\Box$\\
\\
\indent
Now we introduce axiom (AIT3), which we call the Lefschetz type formula.\\
\\ 
{\bf (AIT3)} 
For each $Y \in {\cal Y}$ and all $n\geq 0$,
$$ \mbox{tr}(F_A(Y)^n)=\langle \Phi_A(Y)^n v_{\delta}(Y), v_{\delta}(Y) \rangle_{V_{\Bbb C}}. $$
\\
Here $\mbox{tr}(F_A(Y)^n)$ denotes the trace of $F_A(Y)^n$.\\ 
\\
\noindent
{\bf 4. Standard models of abstract intersection theory}\\
\\
\indent
In this section we construct a model
$$ {\Bbb A}{\Bbb I}{\Bbb T}_m=(V_m,v_{01,m},v_{10,m},v_{\delta,m},h_{a,m},\beta_m,\Phi_{A,m},F_{A,m},{\cal H}_m) $$
of an abstract intersection theory ${\Bbb A}{\Bbb I}{\Bbb T}$.
We call ${\Bbb A}{\Bbb I}{\Bbb T}_m$ which satisfies (AIT1)--(AIT3) a {\it standard model} of abstract intersection theory.

Recall that we have constructed the models $F_{A,m}$ and ${\cal H}_m$ of $F_A$ and ${\cal H}$ in $\S$2. We will construct the remaining elements of the model below.  
 
Let $\{e_i\}_{i=1}^N$ ($1 \leq N:=\dim_{\Bbb C}H \leq \infty$) be an orthonormal basis of the ${\Bbb C}$-Hilbert space $H$. Therefore
$$ H=\Bigl\{ \sum_{i=1}^N \alpha_i e_i \Big\vert \alpha_i \in {\Bbb C}, \sum_{i=1}^N \vert \alpha_i \vert^2 < \infty \Bigr\}. $$
Let $H^1$ be an ${\Bbb R}$-Hilbert space defined by
$$ H^1:=\Bigl\{ \sum_{i=1}^N \alpha_i e_i \Big\vert \alpha_i \in {\Bbb R}, \sum_{i=1}^N \vert \alpha_i \vert^2 < \infty \Bigr\}. $$
Then we have $H^1_{\Bbb C}(:=H^1 \otimes_{\Bbb R} {\Bbb C})=H$ by identifying $e_i \otimes \alpha$ with $\alpha e_i$ for $\alpha \in {\Bbb C}$.

Define ${\Bbb R}$-linear spaces $H^0$ and $H^2$ by
$$ H^0:=\{ \alpha f\vert \alpha \in {\Bbb R} \}\quad \mbox{and}\quad H^2:=\{ \alpha g\vert \alpha \in {\Bbb R} \} $$
with
$$ \langle f, f \rangle _{H^0}:=0\quad \mbox{and}\quad \langle g, g \rangle_{H^2}:=0. $$
\\
{\it Remark 4.1.} The reason why $f\in H^0$ and $g \in H^2$ are defined so that they satisfy the above conditions for degenerate inner product is that 
(IP-b) and (IP-c) in $\S$3 must be satisfied. See (IP-b) and (IP-c) in the proof of Lemma 4.1 below. \hfill $\Box$\\

Then the complexifications $H^0_{\Bbb C}:=H^0 \otimes_{\Bbb R} {\Bbb C}$ and $H^2_{\Bbb C}:=H^2 \otimes_{\Bbb R} {\Bbb C}$ are regarded naturally as
$$ H^0_{\Bbb C} =\{ \alpha f\vert \alpha \in {\Bbb C} \}\quad \mbox{and}\quad H^2_{\Bbb C} =\{ \alpha g\vert \alpha \in {\Bbb C} \} $$
by identifying $f \otimes \alpha$ (resp.\,$g \otimes \alpha$) with $\alpha f$ (resp.\,$\alpha g$) for $\alpha \in {\Bbb C}$.

Let
$$ H^{\bullet}:=H^0 \oplus H^1 \oplus H^2. $$
Here $\oplus$ means the orthogonal direct sum.
That is, we assume that $f$ and $g$ are linearly independent and that 
$\langle f, x \rangle_{H^{\bullet}}=\langle x, f \rangle_{H^{\bullet}}=0$ for $x \in H^1 \oplus H^2$ and
$\langle g, x \rangle_{H^{\bullet}}=\langle x, g \rangle_{H^{\bullet}}=0$ for $x \in H^0 \oplus H^1$.
The inner product $\langle \cdot, \cdot \rangle_{H^{\bullet}}$ on $H^{\bullet}$ is inherited from $\langle \cdot, \cdot \rangle_{H^i}$ ($i=0,1,2$), that is
$\langle x_i, y_i \rangle_{H^{\bullet}}:=\langle x_i, y_i \rangle_{H^i}$ for $x_i, y_i \in H^i$. 

Define an ${\Bbb R}$-linear space $V_m$ by
$$ V_m:=(H^0 \otimes_{\Bbb R} H^2) \oplus (H^1 \otimes_{\Bbb R} H^1) \oplus (H^2 \otimes_{\Bbb R} H^0) $$
with
$$ \langle x_1\otimes x_2, y_1\otimes y_2 \rangle_{V_m}:=\langle x_1, y_1 \rangle_{H^{\bullet}} \langle x_2, y_2 \rangle_{H^{\bullet}}. $$

Since
$$ H^1 \otimes_{\Bbb R} H^1 
= \Bigl\{ \sum_{i,j=1}^N \alpha_{ij} e_i \otimes e_j \Big\vert \alpha_{ij}\in {\Bbb R}, \sum_{i,j=1}^N \vert \alpha_{ij}\vert^2 <\infty \Bigr\}, $$
we have 
$$ (H^1 \otimes_{\Bbb R} H^1)_{\Bbb C} =
\Bigl\{ \sum_{i,j=1}^N \alpha_{ij} e_i \otimes e_j \Big\vert \alpha_{ij}\in {\Bbb C}, \sum_{i,j=1}^N \vert \alpha_{ij}\vert^2 <\infty \Bigr\}
=H^1_{\Bbb C} \otimes_{\Bbb C} H^1_{\Bbb C} $$ 
by identifying $(e_i \otimes e_j) \otimes \alpha$ with $\alpha e_i \otimes e_j$ for $\alpha \in {\Bbb C}$. 
Note that $\{ e_i \otimes e_j \}_{i,j=1}^N$ is an orthonormal basis of the tensor products $H^1 \otimes_{\Bbb R} H^1$ 
and $H^1_{\Bbb C} \otimes_{\Bbb C} H^1_{\Bbb C}$. 
Similarly, by identifying $(f \otimes g) \otimes \alpha$ (resp.\,$(g \otimes f) \otimes \alpha$) with $\alpha f \otimes g$ (resp.\,$\alpha g \otimes f$)
for $\alpha \in {\Bbb C}$, we have  
$$ (H^0 \otimes_{\Bbb R} H^2)_{\Bbb C} = \{ \alpha f \otimes g \vert \alpha \in {\Bbb C} \} = H^0_{\Bbb C} \otimes_{\Bbb C} H^2_{\Bbb C} $$
and
$$ (H^2 \otimes_{\Bbb R} H^0)_{\Bbb C} = \{ \alpha g \otimes f \vert \alpha \in {\Bbb C} \} = H^2_{\Bbb C} \otimes_{\Bbb C} H^0_{\Bbb C}. $$
Note that generally we have $(X \otimes_{\Bbb R} Y)_{\Bbb C} = X_{\Bbb C} \otimes_{\Bbb C} Y_{\Bbb C}$.

Therefore we now have
$$ (V_m)_{\Bbb C} = V_m \otimes_{\Bbb R} {\Bbb C} =
(H^0_{\Bbb C} \otimes_{\Bbb C} H^2_{\Bbb C}) \oplus (H^1_{\Bbb C} \otimes_{\Bbb C} H^1_{\Bbb C}) \oplus (H^2_{\Bbb C} \otimes_{\Bbb C} H^0_{\Bbb C}) $$ 
with
$$ \langle x_1\otimes x_2, y_1\otimes y_2 \rangle_{(V_m)_{\Bbb C}}
=\langle x_1, y_1 \rangle_{H^{\bullet}_{\Bbb C}} \langle x_2, y_2 \rangle_{H^{\bullet}_{\Bbb C}}, $$
where
$$ H^{\bullet}_{\Bbb C}=H^0_{\Bbb C} \oplus H^1_{\Bbb C} \oplus H^2_{\Bbb C} $$
as the orthogonal direct sum. 
Note that the complexification $\langle \cdot, \cdot \rangle_{H^{\bullet}_{\Bbb C}}$ of the inner product $\langle \cdot, \cdot \rangle_{H^{\bullet}}$ is given by
$\langle \alpha_1 x_1, \alpha_2 x_2 \rangle_{H^{\bullet}_{\Bbb C}}:=\alpha_1 \overline{\alpha_1} \langle x_1, x_2 \rangle_{H^{\bullet}}$ for
$x_1, x_2 \in H^{\bullet}$ and $\alpha_1, \alpha_2 \in {\Bbb C}$.

Extend the operator $A$ on $H^1_{\Bbb C}(=H)$ to the operator $A$ on $H^{\bullet}_{\Bbb C}$ by
$$ Af=A\vert_{H^0_{\Bbb C}}f:=0 \quad \mbox{and} \quad Ag=A\vert_{H^2_{\Bbb C}}g:=g. $$

Accordingly, we extend the map $F_{A,m}\colon {\cal Y} \to B(H)$ 
to $F_{A,m}\colon {\cal Y} \to \mbox{End}_{\Bbb C}(H^{\bullet}_{\Bbb C})$ so that
$$ F_{A,m}(Y)f:=e^{t(Y)A\vert_{H^0_{\Bbb C}}}f=f~~ \mbox{and}~~ F_{A,m}(Y)g:=e^{t(Y)A\vert_{H^2_{\Bbb C}}}g=q(Y)g $$
for $Y \in {\cal Y}$. Here $\mbox{End}_{\Bbb C}(H^{\bullet}_{\Bbb C})$ denotes the set of ${\Bbb C}$-linear operators $T\colon H^{\bullet}_{\Bbb C} \supset 
\mbox{dom}(T) \to H^{\bullet}_{\Bbb C}$ with $\mbox{dom}(T)=H^{\bullet}_{\Bbb C}$.\\
\\ 
Let 
$$ v_{01,m}:=f \otimes g \in H^0 \otimes_{\Bbb R} H^2 \subset H^0_{\Bbb C} \otimes_{\Bbb C} H^2_{\Bbb C} \quad \mbox{and} \quad 
v_{10,m}:=g \otimes f \in H^2 \otimes_{\Bbb R} H^0 \subset H^2_{\Bbb C} \otimes_{\Bbb C} H^0_{\Bbb C}. $$

Recall that ${\cal H}_m(Y):=\mbox{Image}(P_{\sigma_Y(A)}) \subset H^1_{\Bbb C}$.
Recall also that $F_{A,m}(Y) {\cal H}_m(Y) \subset {\cal H}_m(Y)$ (i.e.,\,(Frob-a)) by Lemma 2.1 (ii). 
Recall that, by (OP2) and (OP3-a), ${\cal H}_m(Y)$ is finite-dimensional. Let $g(Y):=\frac{1}{2}\dim_{\Bbb C} {\cal H}_m(Y)$. 
Let $\{e_i^{Y}\}_{i=1}^{2g(Y)}$ be an orthonormal basis of ${\cal H}_m(Y)$.

For each $Y \in {\cal Y}$ let 
$$ v_{\delta,m}(Y):=\Bigl(\sum_{i=1}^{2g(Y)} e^{Y}_i \otimes e^{Y}_i\Bigr) +v_{01,m} +v_{10,m} \in (V_m)_{\Bbb C}. $$

Let $\Phi_{A,m}(Y):=I \otimes F_{A,m}(Y)$, where $I$ denotes the identity operator on $H^{\bullet}_{\Bbb C}=H^0_{\Bbb C} \oplus H^1_{\Bbb C} \oplus H^2_{\Bbb C}$.\\
\\
\noindent
{\bf Lemma 4.1.} 
{\it Suppose that an operator} $A\colon H \supset \mbox{dom}(A) \to H$ 
{\it that satisfies} (OP1), (OP2) {\it and} (OP3-a) 
{\it is given. Then the above construction satisfies\\} 
(i) {\it The conditions} (IP-a)--(IP-f).\\
(ii) {\it The Lefschetz type formula} (AIT3).\\ 
\\
\noindent
{\it Proof.} (i) (IP-a) is obvious from definition.\\ 
(IP-b): $\langle v_{01,m}, v_{01,m} \rangle_{V_m} = \langle f \otimes g, f \otimes g \rangle_{V_m} 
= \langle f, f \rangle_{H^{\bullet}} \langle g, g \rangle_{H^{\bullet}}=0$.\\
(IP-c): $\langle v_{10,m}, v_{10,m} \rangle_{V_m} = \langle g \otimes f, g \otimes f \rangle_{V_m} 
= \langle g, g \rangle_{H^{\bullet}} \langle f, f \rangle_{H^{\bullet}}=0$.\\
(IP-d): $\langle v_{01,m}, v_{10,m} \rangle_{V_m} = \langle f \otimes g, g \otimes f \rangle_{V_m} 
= \langle f, g \rangle_{H^{\bullet}} \langle g, f \rangle_{H^{\bullet}}=0$.\\
Since $F_{A,m}(Y)^nf=f$, $F_{A,m}(Y)^ng=q(Y)^ng$ and ${\cal H}_m(Y)$ is $F_{A,m}(Y)$-invariant, we have
\begin{eqnarray*}
\Phi_{A,m}(Y)^n v_{\delta,m}(Y) &=& I \otimes F_{A,m}(Y)^n \{ \sum_{i=1}^{2g(Y)} e_i^Y \otimes e_i^Y + f \otimes g + g \otimes f \}\\
                        &=& \sum_{i=1}^{2g(Y)} e_i^Y \otimes F_{A,m}(Y)^n e_i^Y + f \otimes F_{A,m}(Y)^n g + g \otimes F_{A,m}(Y)^n f\\
                        &=& \sum_{i=1}^{2g(Y)} e_i^Y \otimes F_{A,m}(Y)^n e_i^Y + f \otimes q(Y)^n g + g \otimes f\\
                        &=& \sum_{i=1}^{2g(Y)} e_i^Y \otimes F_{A,m}(Y)^n e_i^Y + q(Y)^n f \otimes g + g \otimes f.
\end{eqnarray*}
(IP-e) and (IP-f) follow from this since $H^0 \perp H^1 \perp H^2$ and $\langle f, f \rangle_{H^0}=\langle g, g \rangle_{H^2}=0$.\\
(ii) To show (AIT3) note that 
\begin{eqnarray*}
 & & \langle \Phi_{A,m}(Y)^n v_{\delta,m}(Y), v_{\delta,m}(Y) \rangle_{(V_m)_{\Bbb C}} \\                                            
 &=& \langle \sum_{i=1}^{2g(Y)} e_i^Y \otimes F_{A,m}(Y)^n e_i^Y + q(Y)^n f \otimes g + g \otimes f,
\sum_{j=1}^{2g(Y)} e_j^Y \otimes e_j^Y + f \otimes g + g \otimes f \rangle_{(V_m)_{\Bbb C}} \\
 &=& \langle \sum_{i=1}^{2g(Y)} e_i^Y \otimes F_{A,m}(Y)^n e_i^Y, \sum_{j=1}^{2g(Y)} e_j^Y \otimes e_j^Y \rangle_{H^1_{\Bbb C} \otimes_{\Bbb C} H^1_{\Bbb C}} \\
 &=& \sum_{i=1}^{2g(Y)} \sum_{j=1}^{2g(Y)} \langle e_i^Y, e_j^Y \rangle_{H^1_{\Bbb C}} \langle F_{A,m}(Y)^n e_i^Y, e_j^Y \rangle_{H^1_{\Bbb C}} \\
 &=& \sum_{i=1}^{2g(Y)} \langle F_{A,m}(Y)^n e_i^Y, e_i^Y \rangle_{H^1_{\Bbb C}} \\
 &=& \mbox{tr}(F_{A,m}(Y)^n). \\
\end{eqnarray*}
This completes the proof of (AIT3). \hfill $\Box$\\
\\
\noindent
{\bf Lemma 4.2.} 
{\it In the same situation as in Lemma 4.1 and its proof, suppose that} (IP-g) {\it further holds.}
{\it Then there is a bilinear form $\beta_m\colon V_m \times V_m \to {\Bbb R}$ and a Hodge vector $h_{a,m}\in V$ which satisfy} 
(AIT1), (AIT2), ($\ast$) {\it and} ($\ast\ast$).\\
\\
\noindent 
{\it Proof.}\\ 
Proof of (AIT1), ($\ast$) and ($\ast\ast$): Recall that
$$ V_m=(H^0 \otimes_{\Bbb R} H^2) \oplus (H^1 \otimes_{\Bbb R} H^1) \oplus (H^2 \otimes_{\Bbb R} H^0), $$  
$$ v_{01,m}=f \otimes g \in H^0 \otimes_{\Bbb R} H^2 \quad \mbox{and} \quad v_{10,m}=g \otimes f \in H^2 \otimes_{\Bbb R} H^0. $$
Therefore we can set
$$ \beta_m(v_{01,m},v_{01,m}):=0, \quad \beta_m(v_{10,m},v_{10,m}):=0 \quad \mbox{and} \quad \beta_m(v_{01,m},v_{10,m})=\beta_m(v_{10,m},v_{01,m}):=1, $$
which are (AIT1-b), (AIT1-c) and (AIT1-d), respectively. Furthermore we can set
$$ \beta_m(x,v_{01,m})=\beta_m(v_{01,m},x):=0 \quad \mbox{and} \quad \beta_m(x,v_{10,m})=\beta_m(v_{10,m},x):=0 $$
for all $x \in H^1 \otimes_{\Bbb R} H^1$. Therefore we have
$$ (\beta_m)_{\Bbb C}(x,v_{01,m})=(\beta_m)_{\Bbb C}(v_{01,m},x)=0 \quad \mbox{and} \quad (\beta_m)_{\Bbb C}(x,v_{10,m})=(\beta_m)_{\Bbb C}(v_{10,m},x)=0 $$
for all $x \in (H^1 \otimes_{\Bbb R} H^1)_{\Bbb C}=H^1_{\Bbb C}\otimes_{\Bbb C}H^1_{\Bbb C}$.

Now for each $Y \in {\cal Y}$ let
$$ v_{\delta 1,m}(Y):=\sum_{i=1}^{2g(Y)} e^{Y}_i \otimes e^{Y}_i \in (H^1 \otimes_{\Bbb R} H^1)_{\Bbb C}=H^1_{\Bbb C}\otimes_{\Bbb C}H^1_{\Bbb C}. $$ 
Note that $v_{\delta,m}(Y)=v_{\delta 1,m}(Y)+v_{01,m}+v_{10,m}$.
Recall from the proof of Lemma 4.1 that 
\begin{eqnarray*}
\Phi_{A,m}(Y)^n v_{\delta,m}(Y) &=& \sum_{i=1}^{2g(Y)} e_i^Y \otimes F_{A,m}(Y)^n e_i^Y + q(Y)^n f \otimes g + g \otimes f \\
                          &=& \sum_{i=1}^{2g(Y)} e_i^Y \otimes F_{A,m}(Y)^n e_i^Y + q(Y)^n v_{01,m} + v_{10,m} \\
                          &=& \Phi_{A,m}(Y)^n v_{\delta 1,m}(Y) + q(Y)^n v_{01,m} + v_{10,m}.
\end{eqnarray*}
Thus, since $\Phi_{A,m}(Y)^n v_{\delta 1,m}(Y) \in (H^1 \otimes_{\Bbb R} H^1)_{\Bbb C}$, we have
$$ (\beta_m)_{\Bbb C}(\Phi_{A,m}(Y)^n v_{\delta,m}(Y), v_{01,m})=1 
\quad \mbox{and} \quad (\beta_m)_{\Bbb C}(\Phi_{A,m}(Y)^n v_{\delta,m}(Y), v_{10,m})=q(Y)^n=O(q(Y)^n), $$
which are (AIT1-e) and (AIT1-f), respectively.
Now that we are given $\beta_m(x,v_{01,m})$, $\beta_m(x,v_{10,m})$, $\beta_m(v_{10,m},y)$ and $\beta_m(v_{01,m},y)$, and $\langle x, y \rangle_{V_m}$,
we can define $\beta_m(x,y)$ for $x, y \in V_m$ by  
$$ \beta_m(x,y):=\beta_m(x,v_{01,m})\beta_m(v_{10,m},y)+\beta_m(x,v_{10,m})\beta_m(v_{01,m},y)-\langle x, y \rangle_{V_m}. $$
Then we see that (AIT1-a), ($\ast$) and ($\ast\ast$) are satisfied. Now we have
$$ (\beta_m)_{\Bbb C}(\Phi_{A,m}(Y)^n v_{\delta,m}(Y),\Phi_{A,m}(Y)^n v_{\delta,m}(Y)) $$
$$ =(\beta_m)_{\Bbb C}(\Phi_{A,m}(Y)^n v_{\delta,m}(Y),v_{01,m})\,\cdot\,(\beta_m)_{\Bbb C}(v_{10,m},\Phi_{A,m}(Y)^n v_{\delta,m}(Y)) $$
$$ +(\beta_m)_{\Bbb C}(\Phi_{A,m}(Y)^n v_{\delta,m}(Y),v_{10,m})\,\cdot\,(\beta_m)_{\Bbb C}(v_{01,m},\Phi_{A,m}(Y)^n v_{\delta,m}(Y)) $$
$$ -\langle \Phi_{A,m}(Y)^n v_{\delta,m}(Y), \Phi_{A,m}(Y)^n v_{\delta,m}(Y) \rangle_{(V_m)_{\Bbb C}}. $$
(AIT1-g) follows from this and (IP-g).\\
Proof of (AIT2): Let $h_{a,m}:=v_{01,m}+v_{10,m}$. If $\beta_m(x,h_{a,m})=0$, then $\beta_m(x,v_{10,m})=-\beta_m(x,v_{01,m})$. 
Hence we have 
$$ \beta_m(x,x)=2\beta_m(x,v_{01,m})\beta_m(x,v_{10,m})-\langle x, x \rangle_{V_m}=-2\beta_m(x,v_{01,m})^2-\langle x, x \rangle_{V_m} \leq 0. $$
Therefore $h_{a,m}$ is a Hodge vector. \hfill $\Box$\\
\\
{\it Remark 4.2.} Note that, given an inner product $\langle \cdot, \cdot \rangle_{H^1_{\Bbb C}}$ for $H^1_{\Bbb C}=H$, the choice of $\beta_m$ is not unique
in our construction of standard models. \hfill $\Box$\\
\\
{\bf 5. Main theorems}\\
\\
\indent
We use the following lemma (see, e.g.,\,[Mon, Lemma 2.2, p.\,20]) in the proof of Theorem 5.2 below.\\
\\
\noindent
{\bf Lemma 5.1.} {\it 
Let $\lambda_i$ $(1\leq i \leq N <\infty)$ be complex numbers. Then there exist infinitely many integers $n\geq 1$ such that $\vert \lambda_1 \vert^n
\leq \vert \sum_{i=1}^{N} \lambda_i^n \vert$.}\\
\\
\noindent
{\bf Theorem 5.2.} {\it Let} $A\colon H \supset \mbox{dom}(A)\to H$ {\it be an operator satisfying} (OP1), (OP2), (OP3-a), (OP4) {\it and} (OP5).\\
(1) {\it If there exists an abstract intersection theory ${\Bbb A}{\Bbb I}{\Bbb T}$ $($in the sense of $\S3$$)$ for $A$, then the Riemann hypothesis holds for $A$.}\\
(2) {\it Suppose further that $A$ satisfies} (OP3-b). {\it Then,
there exists a standard model ${\Bbb A}{\Bbb I}{\Bbb T}_m$ for $A$ if and only if the Riemann hypothesis holds for
$A$ and $A$ is semi-simple.}\\
\\
{\it Proof.} 
(1): Suppose that the RH for $A$ does not hold. Then by (OP5) one can find and fix $Y \in {\cal Y}$ so that $\sigma_Y(A)$ contains 
$s_{\alpha}, s_{\beta} \in \sigma(A)$ with $\mbox{Re}(s_{\alpha})<\frac12, \mbox{Re}(s_{\beta})>\frac12$, respectively. 
Therefore $\sigma_Y(A)$ contains $s_1$ such that $q(Y)^{{\tiny \mbox{Re}}(s_1)}>q(Y)^{\frac12}$.  
Actually, if $0<q(Y)<1$ set $s_1=s_{\alpha}$, while if $q(Y)>1$ set $s_1=s_{\beta}$.

Recall that $\sigma_Y(A)$ is a finite set.
Let $s_i\,(2\leq i \leq 2g(Y):=\dim_{\Bbb C}{\cal H}(Y))$ be all the other eigenvalues of $A$ in $\sigma_Y(A)$, counted with algebraic multiplicities. 
Let $\lambda_i=q(Y)^{s_i}$\,($1\leq i \leq 2g(Y)$). 
Then by Lemma 5.1, $\nu_n := \sum_{i=1}^{2g(Y)} \lambda_i^n$ is not $O(q(Y)^{\frac{n}{2}})$, since we could choose $s_1$ so that
$\vert \lambda_1 \vert^n=\vert q(Y)^{s_1} \vert^n=q(Y)^{\frac{n}{2}}(1+\epsilon)^n$ for some $\epsilon>0$. 

By (Frob-b), we have
$$ \sigma(F_A(Y)^n)=\sigma_p(F_A(Y)^n) = \{q(Y)^{ns} \vert s \in \sigma_Y(A)\} \cup \{0\} = \{\lambda_i^n\vert 1\leq i \leq 2g(Y)\} \cup \{0\}. $$
By (AIT3) and Lemma 3.3 (the Cauchy-Schwarz inequality), we have
\begin{eqnarray*}
\vert \nu_n \vert=\vert\mbox{tr}(F_A(Y)^n)\vert &=& \vert\langle \Phi_A(Y)^n v_{\delta}(Y), v_{\delta}(Y) \rangle_{V_{\Bbb C}}\vert \\
                                                &\leq& \sqrt{\vert\langle v_{\delta}(Y), v_{\delta}(Y) \rangle_{V_{\Bbb C}}\vert \cdot
                                                       \vert\langle \Phi_A(Y)^n v_{\delta}(Y), \Phi_A(Y)^n v_{\delta}(Y) \rangle_{V_{\Bbb C}}\vert}.
\end{eqnarray*}                 
Therefore, by (IP-g), we see that $\nu_n$ is $O(q(Y)^{\frac{n}{2}})$. However, this is a contradiction.\\
\indent
{\it If} part of (2): By Lemma 4.1, we have (IP-a)--(IP-f) and (AIT3) for $V_m$ and $\Phi_{A,m}(Y)$. Therefore all we have to do now is to verify (IP-g)
to apply Lemma 4.2.
Since the RH for the operator $A$ is assumed to hold, each eigenvalue $\lambda_{\ell}$ ($1\leq \ell \leq 2g(Y)$), counted with algebraic multiplicities, of $F_{A,m}(Y)$ can be written as $\lambda_{\ell}=q(Y)^{\frac12}e^{i\theta_{\ell}}$ ($\theta_{\ell} \in {\Bbb R}$).
By the semi-simplicity assumption for $A$, one can choose eigenvectors $w_{\ell}$ associated with
$\lambda_{\ell}$ so that $F_{A,m}(Y) w_{\ell}=\lambda_{\ell} w_{\ell}$. 
Recall that $\{e_i^Y\}_{i=1}^{2g(Y)}$ ($g(Y):=\frac12 \dim_{\Bbb C} {\cal H}_m(Y)$) is an orthonormal basis of ${\cal H}_m(Y)$ (see $\S$4). 
Now one can write $e_i^Y$ as $e_i^Y=\sum_{\ell=1}^{2g(Y)} \alpha_{i \ell}w_{\ell}$ for some $\alpha_{i \ell} \in {\Bbb C}$. Then
\begin{eqnarray*}
 & & \langle \Phi_{A,m}(Y)^n v_{\delta,m}(Y), \Phi_{A,m}(Y)^n v_{\delta,m}(Y) \rangle_{(V_m)_{\Bbb C}} \\                                            
 &=& \langle \sum_{i=1}^{2g(Y)} e_i^Y \otimes F_{A,m}(Y)^n e_i^Y + q(Y)^n f \otimes g + g \otimes f,
\sum_{j=1}^{2g(Y)} e_j^Y \otimes F_{A,m}(Y)^n e_j^Y + q(Y)^n f \otimes g + g \otimes f \rangle_{(V_m)_{\Bbb C}} \\
 &=& \langle \sum_{i=1}^{2g(Y)} e_i^Y \otimes F_{A,m}(Y)^n e_i^Y, 
             \sum_{j=1}^{2g(Y)} e_j^Y \otimes F_{A,m}(Y)^n e_j^Y \rangle_{H^1_{\Bbb C} \otimes_{\Bbb C} H^1_{\Bbb C}} \\
 &=& \sum_{i=1}^{2g(Y)} \sum_{j=1}^{2g(Y)} \langle e_i^Y, e_j^Y \rangle_{H^1_{\Bbb C}} \langle F_{A,m}(Y)^n e_i^Y, F_{A,m}(Y)^n e_j^Y \rangle_{H^1_{\Bbb C}} \\
 &=& \sum_{i=1}^{2g(Y)} \langle F_{A,m}(Y)^n e_i^Y, F_{A,m}(Y)^n e_i^Y \rangle_{H^1_{\Bbb C}} \\
 &=& \sum_{i=1}^{2g(Y)} \langle \sum_{\ell=1}^{2g(Y)} \alpha_{i \ell} F_{A,m}(Y)^n w_{\ell}, 
                                                                                   \sum_{m=1}^{2g(Y)} \alpha_{i m}F_{A,m}(Y)^n w_m \rangle_{H^1_{\Bbb C}}
\end{eqnarray*}
Since $F_{A,m}(Y)^n w_{\ell}=\lambda_{\ell}^n w_{\ell}$, we have (IP-g). 
Therefore by Lemma 4.2, we have (AIT1) and (AIT2) for $V_m$. 

{\it Only if} part of (2): By Lemma 2.1 (i), (ii) and (iii), ${\Bbb A}{\Bbb I}{\Bbb T}_m \Rightarrow \mbox{RH}$ can be proved as (1). 

Let us now show ${\Bbb A}{\Bbb I}{\Bbb T}_m \Rightarrow$ semi-simplicity. Suppose that we have ${\Bbb A}{\Bbb I}{\Bbb T}_m$ but
$A$ is not semi-simple to the contrary. Then one can find and fix $Y \in {\cal Y}$ such that
$$ \sigma_Y(A) = \{ s_1, s_2, \ldots, s_{N-1}, s_N \} $$
which satisfies
$$ \vert \mbox{Im}(s_1) \vert < \vert \mbox{Im}(s_2) \vert < \cdots < \vert \mbox{Im}(s_{N-1}) \vert < \vert \mbox{Im}(s_N) \vert $$
with
$$ \nu(s_1) = \nu(s_2) = \cdots = \nu(s_{N-1}) =1 \quad \mbox{and} \quad \nu(s_N) > 1. $$
What we want to do is to calculate
$$ \langle \Phi_{A,m}(Y)^n v_{\delta,m}(Y), \Phi_{A,m}(Y)^n v_{\delta,m}(Y) \rangle_{(V_m)_{\Bbb C}} $$
and show that it is not of order $O(q(Y)^n)$, which contradicts (IP-g).

We use the notation in Lemma 2.1 and its proof. Then $m_i=\dim_{\Bbb C} H(s_i)=\nu(s_i)=\mbox{mult}(s_i)$ ($1 \leq i \leq N$) by (OP3-b).
We regard $H(s_i)$ as ${\Bbb C}^{m_i}$, that is $H(s_i) \simeq {\Bbb C}^{m_i}$.  
Then we can take a basis $w_{i,\ell} \in {\Bbb C}^{m_i}$ ($1 \leq \ell \leq m_i=\dim_{\Bbb C} H(s_i)$) of the form  
$$ w_{i,\ell}= 
\left(\begin{array}{c}
0\\
\vdots\\
0\\
1\\
0\\
\vdots\\
0
\end{array}\right) \cdots \ell
$$
so that $e^{t(Y)A\vert_{H(s_i)}}$ can be written in the matrix form $N(s_i)$ as in Lemma 2.1 (v). In other words, $w_{i,\ell}$ ($1 \leq \ell \leq m_i$) are generalized eigenvectors of $M(s_i)=A\vert_{H(s_i)}$. 
Let
$$ J_i :=
\left(\begin{array}{ccccc}
0         & 1   &        &        & \bigzerou \\ 
          & 0   & 1      &        &           \\
          &     & \ddots & \ddots &           \\
          &     &        & 0      & 1         \\ 
\bigzerol &     &        &        & 0   
      \end{array}\right) \in M_{m_i}({\Bbb C}). $$
Since $J_i^m=0$ for $m \geq m_i$, $N(s_i)$ in Lemma 2.1 (v) can be written as
\vspace{1.5cm}
$$ \leqno{(5.1)} $$
\vspace{-2.8cm}
\begin{eqnarray*} 
N(s_i) &=& e^{t(Y)s_i}
\left(\begin{array}{ccccc}
\frac{t(Y)^0}{0!} & \frac{t(Y)^1}{1!} & \cdots                       & \cdots                       & \frac{t(Y)^{m_{i}-1}}{(m_{i}-1)!} \\ 
                             & \frac{t(Y)^0}{0!} & \frac{t(Y)^1}{1!} & \cdots                       & 
\frac{t(Y)^{m_{i}-2}}{(m_{i}-2)!} \\
                             &                              & \ddots                       & \ddots                       & \vdots                       \\
                             &                              &                              & \frac{t(Y)^0}{0!} & \frac{t(Y)^1}{1!} \\ 
\bigzerol                    &                              &                              &                              & \frac{t(Y)^0}{0!}
\end{array}\right) \\
 &=& e^{t(Y)s_i} \sum_{k=0}^{m_i-1}\frac{(t(Y)J_i)^k}{k!} \\
 &=& e^{t(Y)s_i} \sum_{k=0}^{\infty}\frac{(t(Y)J_i)^k}{k!} \\
 &=& e^{t(Y)s_i} e^{t(Y)J_i}.
\end{eqnarray*}
Note that 
$$ e^{nt(Y)J_i}w_{i,\ell}=
\left(\begin{array}{c}
\frac{(nt(Y))^{\ell -1}}{(\ell -1)!} \\
\frac{(nt(Y))^{\ell -2}}{(\ell -2)!} \\
\vdots \\
\frac{(nt(Y))^0}{0!} \\
0 \\
\vdots \\
0
\end{array}\right) = \sum_{k=1}^{\ell} \frac{(nt(Y))^{k-1}}{(k-1)!}w_{i,\ell -k+1} \leqno{(5.2)} $$
for $1 \leq i \leq N$ and $1 \leq \ell \leq m_i$. 

Recall that $\{ e_{\mu}^Y \}_{\mu=1}^{2g(Y)}$ is an orthonormal basis of ${\cal H}_m(Y)$.
Thus $e_{\mu}^Y$ can be written as
$$ e_{\mu}^Y = \sum_{i=1}^N \sum_{\ell =1}^{m_i} \alpha_{i,\ell}^{\mu} w_{i,\ell} $$
for some $\alpha_{i,\ell}^{\mu} \in {\Bbb C}$. Then we have by Lemma 2.1 (iv) and (v)
$$ F_{A,m}(Y)^n e_{\mu}^Y = \sum_{i=1}^N F_{A,m}(Y)^n \sum_{\ell =1}^{m_i} \alpha_{i,\ell}^{\mu} w_{i,\ell}
= \sum_{i=1}^N N(s_i)^n \sum_{\ell =1}^{m_i} \alpha_{i,\ell}^{\mu} w_{i,\ell}
= \sum_{i=1}^N \sum_{\ell =1}^{m_i} \alpha_{i,\ell}^{\mu} N(s_i)^n w_{i,\ell}. $$  
Recall from the proof of the {\it If} part of (2) that
$$ \langle \Phi_{A,m}(Y)^n v_{\delta,m}(Y), \Phi_{A,m}(Y)^n v_{\delta,m}(Y) \rangle_{(V_m)_{\Bbb C}}
= \sum_{\mu=1}^{2g(Y)} \langle F_{A,m}(Y)^n e_{\mu}^Y, F_{A,m}(Y)^n e_{\mu}^Y \rangle_{H^1_{\Bbb C}}. $$
Now using (5.1) and (5.2) we have
\begin{eqnarray*}
 & & \langle F_{A,m}(Y)^n e_{\mu}^Y, F_{A,m}(Y)^n e_{\mu}^Y \rangle_{H^1_{\Bbb C}} \\
 &=& \langle \sum_{i=1}^N \sum_{\ell =1}^{m_i} \alpha_{i,\ell}^{\mu} N(s_i)^n w_{i,\ell}, 
     \sum_{j=1}^N \sum_{m =1}^{m_i} \alpha_{j,m}^{\mu} N(s_j)^n w_{j,m} \rangle_{H^1_{\Bbb C}} \\
 &=& \sum_{i=1}^N \sum_{j=1}^N \sum_{\ell=1}^{m_i} \sum_{m=1}^{m_i} \alpha_{i,\ell}^{\mu} \overline{\alpha_{j,m}^{\mu}}
     \langle N(s_i)^n w_{i,\ell}, N(s_j)^n w_{j,m} \rangle_{H^1_{\Bbb C}} \\
 &=& \sum_{i=1}^N \sum_{j=1}^N \sum_{\ell=1}^{m_i} \sum_{m=1}^{m_i} \alpha_{i,\ell}^{\mu} \overline{\alpha_{j,m}^{\mu}}
     \langle e^{nt(Y)s_i}e^{nt(Y)J_i} w_{i,\ell}, e^{nt(Y)s_j}e^{nt(Y)J_j} w_{j,m} \rangle_{H^1_{\Bbb C}} \\
 &=& \sum_{i=1}^N \sum_{j=1}^N \sum_{\ell=1}^{m_i} \sum_{m=1}^{m_i} \alpha_{i,\ell}^{\mu} \overline{\alpha_{j,m}^{\mu}}
     e^{nt(Y)(s_i + \overline{s_j})} \langle e^{nt(Y)J_i} w_{i,\ell}, e^{nt(Y)J_j} w_{j,m} \rangle_{H^1_{\Bbb C}} \\
 &=& \sum_{i=1}^N \sum_{j=1}^N \sum_{\ell=1}^{m_i} \sum_{m=1}^{m_i} \alpha_{i,\ell}^{\mu} \overline{\alpha_{j,m}^{\mu}}
     e^{nt(Y)(s_i + \overline{s_j})} 
     \Big\langle \sum_{a=1}^{\ell} \frac{(nt(Y))^{a-1}}{(a-1)!} w_{i,\ell-a+1}, \sum_{b=1}^m \frac{(nt(Y))^{b-1}}{(b-1)!} w_{j,m-b+1} \Big\rangle_{H^1_{\Bbb C}}.
\end{eqnarray*}
Let $M_{\mu}:=\mbox{max}\,\{\ell \vert \alpha_{N,\ell}^{\mu} \neq 0\}$. Then $\alpha_{N,\ell}^{\mu} \overline{\alpha_{N,m}^{\mu}}=0$ 
if $\ell > M_{\mu}$ or $m > M_{\mu}$. Note that $\mbox{Re}(s_i)=\frac12$ ($\forall i$) since the RH holds by ${\Bbb A}{\Bbb I}{\Bbb T}_m \Rightarrow$ RH.
Recall that $q(Y)=e^{t(Y)}$. Therefore we have
\begin{eqnarray*}
 & & \langle F_{A,m}(Y)^n e_{\mu}^Y, F_{A,m}(Y)^n e_{\mu}^Y \rangle_{H^1_{\Bbb C}} \\
 &=& \alpha_{N,M_{\mu}}^{\mu} \overline{\alpha_{N,M_{\mu}}^{\mu}} e^{2\mbox{\tiny Re}(s_N)nt(Y)}
     \frac{(nt(Y))^{2(M_{\mu}-1)}}{\{(M_{\mu}-1)!\}^2} \langle w_{N,1}, w_{N,1} \rangle_{H^1_{\Bbb C}} + O(e^{nt(Y)}(nt(Y))^{2M_{\mu}-3}) \\
 &=& \vert \alpha_{N,M_{\mu}}^{\mu} \vert^2 q(Y)^n \frac{(n\log q(Y))^{2(M_{\mu}-1)}}{\{(M_{\mu}-1)!\}^2} \Vert w_{N,1} \Vert_{H^1_{\Bbb C}}^2
     + O(q(Y)^n(nt(Y))^{2M_{\mu}-3}) \\
 &=& C_{\mu} q(Y)^n n^{2(M_{\mu}-1)} + O(q(Y)^n n^{2M_{\mu}-3}),
\end{eqnarray*}
where 
$$ C_{\mu} = \vert \alpha_{N,M_{\mu}}^{\mu} \vert^2 \frac{(\log q(Y))^{2(M_{\mu}-1)}}{\{(M_{\mu}-1)!\}^2} \Vert w_{N,1} \Vert_{H^1_{\Bbb C}}^2 > 0. $$
Let $M:=\mbox{max}\,\{M_{\mu}\vert 1 \leq \mu \leq 2g(Y)\}$. 
Since $e_{\mu}^Y$ ($1 \leq \mu \leq 2g(Y)$) are a basis of ${\cal H}_m(Y)$, $\alpha^{\mu}_{N,m_N}\neq 0$ for at least one $\mu$. 
Hence we have $M=m_N >1$. Now we have
$$ \sum_{\mu=1}^{2g(Y)} \langle F_{A,m}(Y)^n e_{\mu}^Y, F_{A,m}(Y)^n e_{\mu}^Y \rangle_{H^1_{\Bbb C}}
= \Bigl( \sum_{M_{\mu}=m_N} C_{\mu} \Bigr) q(Y)^n n^{2(m_N-1)} + O(q(Y)^n n^{2m_N-3}) \neq O(q(Y)^n), $$
which contradicts (IP-g). This completes the proof. \hfill $\Box$\\
\\
\indent
In our previous paper [BU] we constructed a model of abstract intersection theory based on an analogue of the GNS (Gelfand-Naimark-Segal) representation.
Let us call this model which satisfies (INT1)--(INT3) in [BU] a GNS model and denote it as 
${\Bbb A}{\Bbb I}{\Bbb T}_{\mbox{\tiny GNS}}$. The method of the proof of the above theorem also applies to this model. Therefore we have the following theorem.\\
\\
\noindent
{\bf Theorem 5.3.} {\it Let} $A\colon H \supset \mbox{dom}(A)\to H$ {\it be an operator satisfying} (OP1), (OP2), (OP3), (OP4) {\it and} (OP5).
{\it Suppose further that $A$ satisfies} (OP5-b) {\it in} [BU].
{\it Then there exists a GNS model} ${\Bbb A}{\Bbb I}{\Bbb T}_{\mbox{\tiny GNS}}$ {\it for $A$ if and only if 
the Riemann hypothesis holds for $A$ and $A$ is semi-simple.}\\
\\
\indent
We say that $L(s,\chi)$ satisfies the Riemann hypothesis if any nontrivial zero $s_i$ of $L(s,\chi)$ satisfies $\mbox{Re}(s_i)=\frac12$.
We say that a nontrivial zero $s_i$ of $L(s,\chi)$ is simple if it is a zero of $L(s,\chi)$ of order one.

Combining Theorems 5.2 and 5.3 with Theorem 4.1 (iv) of [U] (see Remark 2.1 (4)) we obtain the following theorem.\\ 
\\
\noindent
{\bf Theorem 5.4.} {\it Let} $A\colon H \supset \mbox{dom}(A) \to H$ {\it be an operator constructed in} [U] {\it corresponding to the Dirichlet $L$-function $L(s,\chi)$ associated with a congruence subgroup $\Gamma$ of $SL_2({\Bbb Z})$. Then}\\
(1) {\it $L(s,\chi)$ satisfies the Riemann hypothesis and its all nontrivial zeros are simple if and only if there exists a standard model} 
${\Bbb A}{\Bbb I}{\Bbb T}_m$ {\it for $A$.}\\
(2) {\it $L(s,\chi)$ satisfies the Riemann hypothesis and its all nontrivial zeros are simple if and only if there exists a GNS model} 
${\Bbb A}{\Bbb I}{\Bbb T}_{\mbox{\tiny GNS}}$ {\it for $A$.}\\
\\
{\it Remark 5.1.} In the above theorem if $\Gamma=SL_2({\Bbb Z})$ then the Dirichlet $L$-function $L(s,\chi)$ reduces to the Riemann zeta-function $\zeta(s)$.
\hfill $\Box$\\
\\
{\bf 6. Analogy with the classical theory}\\
\\
\indent
Recall that Weil's explicit formula (according to Patterson [Pat]) reads as follows:
$$ \underbrace{\phi(0) + \phi(1) - \sum_{\rho} \phi(\rho)}_{\mbox{Spectral~term}} 
= \underbrace{W_{\infty}(f) + \sum_{p\colon {\scriptsize \mbox{prime}}} \log p \sum_{n=1}^{\infty} \{ f(p^n)+f(p^{-n}) \} p^{-\frac{n}{2}}}_{\mbox{Geometric~term}}. $$
Here $f$ is a fast decreasing function on ${\Bbb R}_+$,  
$\phi$ is the Mellin transform of $f$, $W_{\infty}$ is an appropriate functional of $f$, and $\rho$ runs over nontrivial zeros of the Riemann zeta-function (or the $L$-function), counted with multiplicity. For the original work of Weil, see [1952b] and [1972] of [W2]. See also [C] and [CM, p.\,344].

The idea of introducing the model $F_{A,m}(Y)$ of an analogue of the Frobenius morphism in this paper is hinted by the spectral side of the above formula. 
By Lemma 2.2 of [BU, p.\,702] there is a function $\phi_Y(s)$ ($Y \in {\cal Y}$) which is analytic in an open set $\Supset \Omega_{\infty}$ such that\\
(i) $\phi_Y(0)=1$,\\
(ii) $\phi_Y(1)=q(Y)$,\\
(iii) $\phi_Y(s_i)=q(Y)^{s_i}$ if $s_i \in \sigma_Y(A)$,\\
\\
(iv)
\vspace{-1.15cm} 
$$ \lim_{s\to s_i} \frac{\phi_Y(s)}{(s-s_i)^{m_i}}=c_{Y,i} \in {\Bbb C}~\mbox{for~some}~c_{Y,i} \neq 0~\mbox{if}~s_i \in \sigma(A)-\sigma_Y(A)~\mbox{with}~\nu(s_i)=m_i. $$

For this $\phi_Y(s)$, let $\phi_Y(A)\colon H\supset \mbox{dom}(\phi_Y(A)) \to H$ be defined by
$$ \phi_Y(A)x:=\lim_{{\scriptstyle T \to \infty} \atop {\scriptstyle T \in {\cal Y}}}
\frac{1}{2\pi i} \Bigl(\oint_{\partial\Omega_T} \phi_Y(s)(sI-A)^{-1}ds\Bigr)x $$
for 
$$ x\in \mbox{dom}(\phi_Y(A)):=\{x\in H \vert \mbox{the~limit}~\phi_Y(A)x~\mbox{exists~in}~H\}. $$
Then it is easy to prove that $\mbox{dom}(\phi_Y(A))=H$ and that
$$ \phi_Y(A)=F_{A,m}(Y). $$
For the proof use $(sI-M(s_i))^{-1}$ in the proof of Lemma 2.1. It is also easy to see that
$$ \mbox{tr}(\phi_Y(A)) = \sum_{s_i \in \sigma_Y(A)} \mbox{mult}(s_i) \phi_Y(s_i). $$

Let $C$ be a smooth projective curve (one-dimensional scheme) over a finite field ${\Bbb F}_q$. 
Let $\mbox{Frob}$ be the Frobenius morphism on $C$. Then $F_A(Y)$ in $\S$2 is an analogue of $\mbox{Frob}$.
 
For $S=C\times C$, the surface over ${\Bbb F}_q$,
let $\Pic(S)$ be its Picard group, which we regard as a ${\Bbb Z}$-module, so as to preserve the analogy with Weil divisors.
The ${\Bbb R}$-linear space $V$ in $\S$3 is modeled on $\Pic(S)\otimes_{\Bbb Z}{\Bbb R}$. The ${\Bbb R}$-bilinear form 
$\beta(\cdot, \cdot)$ in $\S$3 is modeled on the ${\Bbb R}$-tensored intersection pairing $i(\cdot, \cdot)$ on $\Pic(S)\otimes_{\Bbb Z}{\Bbb R}$.

The operator $\Phi_A(Y)$ in (AIT1) is an analogue of the linear map on $\Pic(S)\otimes_{\Bbb Z}{\Bbb R}$ induced by the morphism 
$\mbox{id} \times \mbox{Frob}$ on $S$.   
Then one may regard $v_{01}$, $v_{10}$, $v_\delta(Y)$ and $\Phi_A(Y)^n v_{\delta}(Y)$ in (AIT1) as analogues of cycles $\mbox{pt} \times C$, $C \times \mbox{pt}$, 
$\Delta$ and $\Gamma_{{\tiny \mbox{Frob}}^n}$ in $\Pic(S)$, respectively.
Here $\Delta$ is the diagonal, and $\Gamma_{{\tiny \mbox{Frob}}^n}$ is the graph of $\mbox{Frob}^n$. So here is the dictionary.

\begin{center}
{\large 
\begin{tabular}{|c|c|} \hline
$\mbox{Pic}(S) \otimes_{\Bbb Z} {\Bbb R}$ & $V$ \\ \hline
$\mbox{pt} \times C$ & $v_{01}$     \\
$C \times \mbox{pt}$ & $v_{10}$     \\
$\Delta$             & $v_{\delta}(Y)$ \\
$\Gamma_{{\tiny \mbox{Frob}}^n}$ & $\Phi_A(Y)^n v_{\delta}(Y)$ \\ \hline
\end{tabular}
}
\end{center}

The cycles $\mbox{pt} \times C$, $C \times \mbox{pt}$, $\Delta$ and $\Gamma_{{\tiny \mbox{Frob}}^n}$ have the following properties.\\
\indent \hspace{0.8cm}
(i) $i(\mbox{pt} \times C, \mbox{pt} \times C)=0$.\quad 
(ii) $i(C \times \mbox{pt}, C \times \mbox{pt})=0$.\quad 
(iii) $i(\mbox{pt} \times C, C \times \mbox{pt})=1$.\\
\indent \hspace{0.8cm}
(iv) $i(\Gamma_{{\tiny \mbox{Frob}}^n}, \mbox{pt} \times C)=1$.\quad 
(v) $i(\Gamma_{{\tiny \mbox{Frob}}^n}, C \times \mbox{pt})=q^n$.\quad 
(vi) $i(\Gamma_{{\tiny \mbox{Frob}}^n}, \Gamma_{{\tiny \mbox{Frob}}^n})=q^n$.\\
The axioms of (AIT1) are analogues of these properties.

The Hodge property in (AIT2) comes from the classical Hodge index theorem. A Hodge vector $h_a$ corresponds to an ample hyperplane section of $S$, 
thereby $\beta(\cdot,h_a)$ gives an analogue of the degree function $\deg \otimes_{\Bbb Z} 1\colon \Pic(S)\otimes_{\Bbb Z} {\Bbb R}\to {\Bbb R}$.
Lemma 3.1 is an analogue of the inequality of Castelnuovo-Severi.

The construction of $V_m$ of a standard model in $\S$4 is hinted by the K\"unneth formula for the \'etale cohomology.
The Tate conjecture for $S=C \times C$ and codimension one is equivalent to that the map 
$\mbox{Pic}\,S \otimes {\Bbb Q}_{\ell} \to H^2_{{\tiny \mbox{\'et}}}(S,{\Bbb Q}_{\ell}(1))$ is bijective (Proposition (4.3) of Tate [T2]).
Note that $H^2_{{\tiny \mbox{\'et}}}(S,{\Bbb Q}_{\ell}(1)) = 
H^2_{{\tiny \mbox{\'et}}}(\bar{S},{\Bbb Q}_{\ell}(1))^{\scriptsize \mbox{Gal}(\overline{{\Bbb F}}_q/{\Bbb F}_q)}$, where 
$\bar{S}=S \times_{{\Bbb F}_q} \overline{{\Bbb F}}_q$ (see [T2]).
Tate [T1] himself has proven his conjecture for abelian varieties over finite fields for the case of codimension one.
From this the Tate conjecture follows also for $S=C\times C$ in the codimension one case.
By the K\"unneth formula for $\ell$-adic cohomology we have
\newpage
$$ H^2_{{\tiny \mbox{\'et}}}(\bar{S},{\Bbb Q}_{\ell}) 
\simeq \Bigl(H^0_{{\tiny \mbox{\'et}}}(\bar{C},{\Bbb Q}_{\ell})\otimes H^2_{{\tiny \mbox{\'et}}}(\bar{C},{\Bbb Q}_{\ell})\Bigr) \oplus
\Bigl(H^1_{{\tiny \mbox{\'et}}}(\bar{C},{\Bbb Q}_{\ell}) \otimes H^1_{{\tiny \mbox{\'et}}}(\bar{C},{\Bbb Q}_{\ell})\Bigr) \oplus 
\Bigl(H^2_{{\tiny \mbox{\'et}}}(\bar{C},{\Bbb Q}_{\ell}) \otimes H^0_{{\tiny \mbox{\'et}}}(\bar{C},{\Bbb Q}_{\ell})\Bigr). $$
Here $\bar{C}=C \times_{{\Bbb F}_q} \overline{{\Bbb F}}_q$.
The definition of the ${\Bbb R}$-linear space $V_m$ is modeled on this. For the K\"unneth formula for $\ell$-adic cohomology see Chap.\,6, $\S$8 of Milne [Mil].

For a morphism $\varphi\colon C \to C$, the Lefschetz fixed-point formula for the $\ell$-adic \'etale cohomology group 
$H^i_{{\tiny \mbox{\'et}}}=H^i_{{\tiny \mbox{\'et}}}(\bar{C}, {\Bbb Q}_{\ell})$ is given by
$$ \mbox{tr}(\varphi^{\ast n}\vert_{H^0_{{\tiny \mbox{\'et}}}}) -\mbox{tr}(\varphi^{\ast n}\vert_{H^1_{{\tiny \mbox{\'et}}}})
+\mbox{tr}(\varphi^{\ast n}\vert_{H^2_{{\tiny \mbox{\'et}}}}) = i(\Gamma_{\varphi^n}, \Delta), $$
where $\Gamma_{\varphi^n}$ is the graph of $\varphi^n$. If $\varphi=\mbox{Frob}$, then it turns out that
$$ 
\mbox{tr}(\varphi^{\ast n}\vert_{H^0_{{\tiny \mbox{\'et}}}})=1=i(\Gamma_{\varphi^n}, \mbox{pt} \times C)i(\Delta, C\times \mbox{pt}) $$ 
and 
$$ \mbox{tr}(\varphi^{\ast n}\vert_{H^2_{{\tiny \mbox{\'et}}}})=q^n=i(\Gamma_{\varphi^n}, C \times \mbox{pt})i(\Delta, \mbox{pt}\times C). $$
So the Lefschetz fixed-point formula reads for $\varphi^n=\mbox{Frob}^n$ as\\
\begin{eqnarray*}
\mbox{tr}(\varphi^{\ast n}\vert_{H^1_{{\tiny \mbox{\'et}}}})  
&=& i(\Gamma_{\varphi^n}, \mbox{pt} \times C)i(\Delta, C\times \mbox{pt})+i(\Gamma_{\varphi^n}, C \times \mbox{pt})i(\Delta, \mbox{pt}\times C)
-i(\Gamma_{\varphi^n}, \Delta)\\
&=:& \langle \Gamma_{\varphi^n}, \Delta \rangle_{\mbox{\scriptsize Pic}(S)\otimes_{\Bbb Z}{\Bbb R}}.
\end{eqnarray*}
(AIT3) is modeled on this. Consider the operators $A$ and $F_{A,m}(Y)$ ($Y \in {\cal Y}$) which are extended to 
$H^{\bullet}_{\Bbb C}=H^0_{\Bbb C} \oplus H^1_{\Bbb C} \oplus H^2_{\Bbb C}$ as in $\S$4. Then we have 
$$ \phi_Y(A)f=F_{A,m}(Y)f=f=\phi_Y(0)f~~ \mbox{and}~~ \phi_Y(A)g=F_{A,m}g=q(Y)g=\phi_Y(1)g $$
for $f \in H^0_{\Bbb C}$ and $g \in H^2_{\Bbb C}$. 
The operator $\phi_Y(A)$ acting on $H^i_{\Bbb C}$ is an analogue of $\mbox{Frob}^{\ast}$ acting on $H^i_{{\tiny \mbox{\'et}}}$ ($i=0,1,2$). Since
\begin{eqnarray*}
\Phi_{A,m}(Y)^n v_{\delta,m}(Y) &=& \sum_{i=1}^{2g(Y)} e_i^Y \otimes F_{A,m}(Y)^n e_i^Y + f \otimes F_{A,m}(Y)^n g + g \otimes F_{A,m}(Y)^n f \\
                                &=& \sum_{i=1}^{2g(Y)} e_i^Y \otimes F_{A,m}(Y)^n e_i^Y + \phi_Y(1)^n v_{01,m} + \phi_Y(0)^n v_{10,m}
\end{eqnarray*}
(see the proof of Lemma 4.2), we have by the setting of the proof of Lemma 4.2
$$ \mbox{tr}(\phi_Y(A)^n\vert_{H^0_{\Bbb C}})=\phi_Y(0)^n
=(\beta_m)_{\Bbb C}(\Phi_{A,m}(Y)^n v_{\delta,m}(Y), v_{01,m})\,\cdot\,(\beta_m)_{\Bbb C}(v_{10,m}, v_{\delta,m}(Y)) $$ 
and 
$$ \mbox{tr}(\phi_Y(A)^n\vert_{H^2_{\Bbb C}})=\phi_Y(1)^n
=(\beta_m)_{\Bbb C}(\Phi_{A,m}(Y)^n v_{\delta,m}(Y), v_{10,m})\,\cdot\,(\beta_m)_{\Bbb C}(v_{01,m}, v_{\delta,m}(Y)). $$
Therefore we have by $(\ast\ast)$
$$ \mbox{tr}(\phi_Y(A)^n\vert_{H^0_{\Bbb C}}) - \mbox{tr}(\phi_Y(A)^n\vert_{H^1_{\Bbb C}}) + \mbox{tr}(\phi_Y(A)^n\vert_{H^2_{\Bbb C}})
= (\beta_m)_{\Bbb C}(\Phi_{A,m}(Y)^n v_{\delta,m}(Y), v_{\delta,m}(Y)), $$ 
which is equivalent to (AIT3).\\
\newpage 
\noindent
{\bf References}\\
\\
{\small
$[$BU$]$ G. Banaszak and Y. Uetake, Abstract intersection theory and operators in Hilbert space, {\it Communications in Number Theory and Physics} {\bf 5}
(2011), 699--712.\\
\\
$[$C$]$ A. Connes, Noncommutative geometry and the Riemann zeta function, {\it Mathematics$:$\,Frontiers and Perspectives}, V. Arnold et al. eds., AMS, 
2000, 35--54.\\
\\
$[$CM$]$
A. Connes and M. Marcolli, {\it Noncommutative Geometry, Quantum Fields and Motives}, AMS Colloquium Publications {\bf 55}, 
AMS, Providence; Hindustan Book Agency, New Delhi, 2008.\\
\\
$[$D$]$ C. Deninger, Some analogies between number theory and dynamical systems on
foliated spaces, {\it Proc.\,of the ICM, Berlin, Vol.\,I 
(Doc.\,Math.\,J.\,DMV)}, 1998, 163--186.\\
\\
$[$Ge$]$
I. M. Gelfand, Automorphic functions and the theory of representations, {\it Proc.\,of the ICM, Stockholm}, 1962, 74--85.\\
\\
$[$GoGK$]$
I. Gohberg, S. Goldberg and M. A. Kaashoek,
    {\it Classes of Linear Operators, Vol.\,I},
    Oper.\,Theory Adv.\,Appl. {\bf 49},
    Birkh\"auser, Basel, 1990.\\
\\    
$[$Gro$]$
A. Grothendieck, Sur une note de Mattuck-Tate, {\it J. reine ang. Math.}
{\bf 200} (1958), 208--215.\\
\\
$[$LP$]$
P. D. Lax and R. S. Phillips, {\it Scattering Theory for Automorphic Functions}, Ann.\,of Math. Studies {\bf 87}, Princeton University Press,
Princeton, 1976.\\
\\
$[$Mac$]$ B. D. MacCluer, {\it Elementary Functional Analysis}, GTM {\bf 253}, Springer, 2009.\\
\\
$[$Mil$]$
J. S. Milne, {\it \'Etale Cohomology}, Princeton Mathematical Series {\bf 33}, Princeton University Press, Princeton, N.J., 1980.\\
\\
$[$Mon$]$
P. Monsky,
    {\it P-Adic Analysis and Zeta Functions},
    Lectures in Mathematics {\bf 4}, Kyoto University, 
    Kinokuniya Book-Store, Tokyo, 1970.\\
\\
$[$Pat$]$
S. J. Patterson, {\it An Introduction to the Theory of the Riemann Zeta-Function}, Cambridge University Press, Cambridge, 1988.\\
\\ 
$[$PavF$]$ 
B. S. Pavlov and L. D. Faddeev, Scattering theory and automorphic functions, {\it Zap. Nau{\v c}n. Sem. Leningrad. Otdel. Mat. Inst. Steklov. (LOMI)}
{\bf 27} (1972), 161--193. (Russian); English transl.:\,{\it J. Soviet Math.} {\bf 3} (1975), 522--548.\\
\\
$[$S$]$
J.-P. Serre, Analogues K\"ahl\'eriens de certaines conjectures de Weil, {\it Ann.\,of Math.}\,{\bf 71} (1960), 392--394.\\
\\
$[$T1$]$ J. Tate, Endomorphisms of abelian varieties over finite fields, {\it Invent.\,Math.}\,{\bf 2} (1966), 134--144.\\ 
\\
$[$T2$]$ 
\_\hspace{-0.05cm}\_\hspace{-0.05cm}\_\hspace{-0.05cm}\_\hspace{-0.05cm}\_\,, Conjectures on algebraic cycles in $\ell$-adic cohomology, {\it Proc.\,of Symposia in Pure Mathematics} {\bf 55}, Part I, 1994, 71--83.\\
\\
$[$U$]$
Y. Uetake, 
Spectral scattering theory for automorphic forms, {\it Integral Equations Operator Theory} {\bf 63} (2009), 439--457.\\
\\
$[$W1$]$ 
A. Weil, {\it \OE uvres Scientifiques Collected Papers}, papers, a private communication and books concerning intersection theory:\,[1940b], [1941], [1942], [1946a], [1948a], [1948b] (Vol.\,I), [1954h] (Vol.\,II), Springer Verlag, New York, 1979.\\
\\
$[$W2$]$ 
\_\hspace{-0.05cm}\_\hspace{-0.05cm}\_\hspace{-0.05cm}\_\hspace{-0.05cm}\_\,, {\it \OE uvres Scientifiques Collected Papers}, papers concerning explicit formulas:\,[1952b] (Vol.\,II), [1972] (Vol.\,III), Springer Verlag, New York, 1979.\\
\\ 
}
\\
\noindent
{\it Faculty of Mathematics and Computer Science\\
Adam Mickiewicz University\\
ul. Umultowska 87, 61-614 Pozna\'n\\
Poland\\}
E-mail: {\tt banaszak@amu.edu.pl, uetake@amu.edu.pl}
\end{document}